\documentclass[reqno,11pt]{amsart}
\usepackage{amsfonts,amsmath,amsthm, url}
\usepackage{amssymb}  
\usepackage{mathrsfs}
\usepackage{graphics}
\usepackage{graphicx}
\usepackage{color}
\usepackage{stmaryrd} 

\usepackage[english]{babel}
\usepackage[T1]{fontenc}
\usepackage{graphicx}
\usepackage{mathrsfs}
\usepackage{mathtools}
\usepackage{float}

\usepackage{cite}
\usepackage{enumerate}
\usepackage{babel}
\usepackage[font=small,labelfont=bf]{caption}

   \newtheorem{thm}{Theorem}
   \newtheorem{prop}{Proposition}

\newenvironment{pr}[1]
   {{\noindent \bf Proof.   }}{\hfill $\Box$}

\newcommand{\Z}{\mathbb{Z}}

\newcommand{\R}{\mathbb{R}}

\newcommand{\E}{\mathbb{E}}

\newcommand{\prob}{\mathbb{P}}

\newcommand{\vertiii}[1]{{\left\vert\kern-0.25ex\left\vert\kern-0.25ex\left\vert #1 
    \right\vert\kern-0.25ex\right\vert\kern-0.25ex\right\vert}}

\makeatletter
  
  \@addtoreset{equation}{section}
\makeatother
\begin{document}

\title[2d stochastic anisotropic Swift-Hohenberg equation]
{Pattern formation in 2d stochastic anisotropic Swift-Hohenberg equation}\textcolor{black}{\footnote{\textcolor{black}{Reika Fukuizumi is grateful to Dao-Zhi Zeng, who served as the editor for this submission.}}}
\author[Reika FUKUIZUMI]{Reika FUKUIZUMI$^{\scriptsize 1}$}
\author[Yueyuan Gao]{Yueyuan Gao$^{\scriptsize 2}$}
\author[Guido Schneider]{Guido Schneider$^{\scriptsize 3}$}
\author[Motomitsu Takahashi]{Motomitsu Takahashi$^{\scriptsize 4}$}



\maketitle

\begin{center} \small
$^{1,4}$ Research Center for Pure and Applied Mathematics, \\
Graduate School of Information Sciences, Tohoku University,\\
Sendai 980-8579, Japan; \\
\end{center}

\begin{center} \small
$^2$ Laboratory of Mathematical Modeling, Research Institute for Electronic Science,\\
 {Hokkaido University}, Sapporo, 060-0812, Japan; \\
\end{center}

\begin{center} \small
$^3$ Institut f\"{u}r Analysis, Dynamik und Modellierung,
Universit\"{a}t Stuttgart, \\
Pfaffenwaldring 57, D-70569 Stuttgart, Germany\\
\end{center}

\vskip 0.1 in
\noindent
{\bf Abstract}. In this paper, we study a phenomenological model for pattern formation in electroconvection, and the effect of noise on the pattern.  As such model we consider an anistropic Swift-Hohenberg equation adding an additive noise. We prove the existence of  a global solution of that equation on the two dimensional torus. In addition, inserting a scaling parameter, we consider the equation on a large domain near
its change of stability. We observe numerically that, under the appropriate scaling, its solutions can be
approximated by a periodic wave, which is modulated by the solutions to a stochastic
Ginzburg-Landau equation.
 
\vskip 0.1 in

\section{Introduction}
The Swift-Hohenberg equation is a celebrated toy model for
the convective instability in the Rayleigh-B\'{e}nard convection \cite{SH}.  This equation has played an important role not only in the model of pattern formation in thermal convection, but also in the study of different fields  including electroconvection, economics, biology, sociology, optics, etc. (see \cite{JK}).
\vspace{3mm}

The one-dimensional Swift-Hohenberg equation is given by
\begin{equation}\label{u_1_dimension}
\partial_tu=-(1+\partial_x^2)^2u+\alpha u-u^3, \,\,\, t\ge 0,\,\,\,  x\in\mathbb{R},
\end{equation}
where $\alpha \in \R$ is called the stress parameter.  
The linear part is clearly analyzed using Fourier transform. The ansatz 
\begin{align*}
u(t,x)=e^{\lambda(k)t+ikx},
\end{align*}
where $k\in\R$ is the wave number, yields $\lambda(k)=-(1-k^2)^2+\alpha$.
If $\alpha>0$, then unstable modes around $k = \pm 1$ exist 
and thus, in this case the convection and the pattern formation occur. 
Now let $\alpha=\varepsilon^2>0$. We expect that the solution can be described by the ansatz 
$$u(t,x)=\varepsilon A(T,X) e^{ix}+c.c.,\, X=\varepsilon x,\, T=\varepsilon^2 t,$$ 
\textcolor{black}{where c.c. means the complex conjugate.}
Substituting it into the above equation, and comparing the coefficients of $\varepsilon^j e^{ikx}(j,k\in\mathbb{Z})$, we see that the so-called residual is minimized if $A(T, X)$ fulfills
\begin{equation}\label{A_1_dimension}
\partial_TA=4\partial_X^2A+A-3|A|^2A.
\end{equation}
Indeed, it is known that if $\varepsilon >0$ is taken to be small enough, $u(t,x)-(A(T,X)e^{ix}+c.c.)$ becomes smaller in a suitable sense (see \cite{HU}).
\vspace{3mm}

We are interested in adding noise in this formulation. 
For the stochastic equation, it is shown in \cite{BBS} 
that the solution $u$ of the one dimensional stochastic Swift-Hohenberg equation; 
\begin{equation}
\partial_tu=-(1+\partial_x^2)^2+\nu\varepsilon^2u-u^3+ \varepsilon^{\frac{3}{2}} \dot{\xi}_{\varepsilon}, 
 \,\,\, t\ge 0,\,\,\,  
\end{equation}
where $\nu>0$, and $\dot{\xi}_{\varepsilon}$ is the real valued space-time white noise,
can be approximated by using the solution $A$ of the stochastic Ginzburg-Landau equation:
\begin{align*}
\partial_TA=4\partial_X^2A+\nu A-3A|A|^2+ \dot{\eta}, ,\, X=\varepsilon x,\, T=\varepsilon^2 t
\end{align*}
where $\dot{\eta}$ is a complex valued noise. The approximation is given by
\begin{equation*}
u(t,x)\simeq\varepsilon A(\varepsilon^2t,\varepsilon x)e^{ix}+c.c.
\end{equation*}
This result is proved on the whole space $\mathbb{R}$,  there is also a result by \cite{BHP}  on the one-dimensional torus. 
To our best knowledge, this approximation in the stochastic case 
is known only in one dimension, the results in more than two dimensions are not known. The main problem in the stochastic case in the dimension more than two is 
that the solution of the stochastic Ginzburg-Landau equation has a priori a negative regularity, thus we need to use a renormalization 
to make sense to the nonlinearity (see for example \cite{MW}), whereas the Swift-Hoheberg equation has good regularity and we do not need to use a renormalization. Rather,  
the linear part of Swift-Hohenberg equation can define the Wick products of Ornstein-Uhlenbeck process in a certain scale of $\varepsilon$. We would address this issue in the sequel paper. 
In the deterministic case, such approximation for more general forms of equation is already known in two dimensions, see for e.g. \cite{GS}. 
\vspace{3mm}

In this paper we thus consider a two-dimensional stochastic Swift-Hohenberg equation on the torus given by
\begin{equation}\label{SSHE}
\begin{cases}
\partial_tu=
-(1+\partial_x^2)^2u+\partial_y^2u+u-u^3+\dot{\xi}, \,\,\, & (x,y) \in \mathbb{T}^2,\\
u(0)=u_0.\,\,\,& 
\end{cases}
\end{equation}
where $\dot{\xi}$ is the real-valued space-time white noise. 
This equation is a phenomenological model for pattern formation in electroconvection, in the sense that 
the spectral surface is similar to the modeling of the electroconvection proposed in \cite{SU}. 
In this paper, as the first step, we prove the existence of a solution of this equation $(\ref{SSHE})$. 
\vspace{3mm}

By the result of the approximation of the deterministic equation and the approximation result of  \cite{BBS}, 
we can expect that the two-dimensional stochastic Swift-Hohenberg equation (\ref{SSHE}) also can be approximated 
by the two-dimensional complex stochastic Ginzburg-Landau equation. We may see by formal 
computations that the solutions of 
$$\partial_tu=
-(1+\partial_x^2)^2u+\partial_y^2u+ \varepsilon^2 u-u^3+ \varepsilon \dot{\xi}_{\varepsilon} $$ 
defined on the domain $[-L/\varepsilon, L/\varepsilon] \times [-L/\varepsilon, L/\varepsilon]$ ($L>0$) 
with periodic boundary condition would be approximated by the solution $A$ of 
$$ \partial_T A =4\partial^2_X A +\partial_Y^2 A +A - 3\lvert A\rvert^2A+\dot{\eta}, $$
on the domain $[-L, L] \times [-L, L]$, through the ansatz 
\begin{equation*}
u(t,x,y) = \varepsilon A(\varepsilon^2t,\varepsilon x, \varepsilon y)e^{ix}+c.c.\end{equation*}
\textcolor{black}{We try to see whether this would be observed or not by numerical simulations.}


The organization of this paper is as follows. 
In Section 2, we prepare the notation necessary for discussing the subsequent sections, and we state our main theorem. In Section 3, 
we investigate the regularity of the solution of the linear equation of  the two-dimensional stochastic Swift-Hohenberg equation using the Kolmogorov test. 
Section 4  is dedicated to prove the existence of the solutions of equation (\ref{SSHE}) using the regularity of the solutions of the linear equation obtained in Section 3. 
We use the compactness method and obtain the solution as the limit of  finite dimensional Galerkin approximation and its energy uniform estimates. 
Lastly, we will present the numerical simulations in Section 5. 

%
%

\section{Preliminaries and main results}

\subsection{Notation}

In this section, we define the notation for our discussion.
\begin{itemize}
\item[(i)] Our results will concern the periodic functions in $\R^2$ and for a fixed $L>0$ we shall take 
the fundamental period in each variable to be $2L$. That is, a function $f$ on $\R^2$ is said to be 
periodic if $f( \mathbf{x}+2L\mathbf{k})=f(\mathbf{x})$ for all $\mathbf{x}=(x,y) \in \R^2$ and $\mathbf{k}=(k,l) \in \Z^2.$ 
For the analysis, a natural option would be to base the definition
of the Sobolev spaces on discrete Fourier series, and those are adapted to the ``torus,'' namely we regard the periodic functions as functions on the space ${\R}^2/((2L)\Z)^2$ which we will call the torus and denote by $\mathbb{T}^2$. 
We identify $\mathbb{T}^2$ with the cube $[-L,L)^2.$ 
\item[(ii)] Let $\nu_2= \frac{1}{(2L)^2} m_2,$ where $m_2$ is two-dimensional Lebesgue measure. 
Then, by the identification above,  $\nu_2$ induces a measure on $\mathbb{T}^2$, but 
we denote it by the same $\nu_2$ with an abuse of notation.
For all $p\in [1,\infty],$ $L^p(\mathbb{T}^2)$ denotes thus $L^p([-L,L]^2)$ with this Lebesgue measure $\nu_2$.
\item[(iii)]  For measurable complex-valued functions $f,g \in L^2(\mathbb{T}^2)$, the $L^2(\mathbb{T}^2)$ inner-product is denoted by 
$$ (f,g) := \int_{\mathbb{T}^2} f(x) \overline{g(x)} d\nu_2= \int_{[-L,L]^2}  f(x) \overline{g(x)} d\nu_2.$$
\item[(iv)] For $\varepsilon \in [0,1]$, we set $\mathcal{L_{\varepsilon}} =-(1+\partial_x^2)^2
+\partial_y^2+\varepsilon^2$ and 
$\lambda_{k,l,\varepsilon} =-\left\{(1- \textcolor{black}{\left( \frac {\pi}{L} \right)^2}  k^2)^2
+ \textcolor{black}{\left( \frac {\pi}{L} \right)^2}  l^2-\varepsilon^2 \right\}$ for $k,l\in \Z$, which are eigenvalues of $\mathcal{L}_{\varepsilon}$. The dependence of the operator $\mathcal{L_{\varepsilon}}$ on  $\varepsilon$ is not essential for the existence of solutions, thus for the sake of simplicity 
we set $\varepsilon=1$, but we use $\varepsilon$ for the purpose of numerical simulations in Section 5.
\item[(v)] Let $\ \{e_{k,l}(x,y)\}_{k,l\in\mathbb{Z}}$ be the eigenfunctions corresponding to $\lambda_{k,l,0}$, which will simply be denoted 
by $\lambda_{k,l}$, i.e., 
$$\mathcal{L}_{0 } e_{k,l} = \lambda_{k,l} e_{k,l},  \quad e_{k,l}(x,y) 
= \textcolor{black} {\frac{1}{2L}} e^{\frac{-i \pi (k,l)}{L} \cdot (x,y)} $$ 
and which constitute a complete orthogonal basis in $L^2(\mathbb{T}^2)$.
\item[(vi)] For $s \in \R$ and $1 \le p <+\infty$, we denote by ${\mathcal{W}}^{s,p}(\mathbb{T}^2)$ the space of  $f \in \mathcal{S}'$ satisfying 
\begin{align*}
\|f\|_{\mathcal{W}^{s,p}(\mathbb{T}^2)} := \|(1-\mathcal{L}_0)^{\frac{s}{2}} f \|_{L^p (\mathbb{T}^2)}.
\end{align*}
\item[(vii)] For $1 \le  p<+\infty$, and $T>0$, and a Banach space $B$ with the norm $||\cdot||_{B}$, we denote by $L^p(0,T;B)$ the $B$-valued measurable function $g$ on $[0,T]$ such that 
\begin{align*}
\int_0^T \|g(t) \|_{B}^p dt<\infty.
\end{align*}
For $\alpha\in(0,1)$, we denote by $W^{\alpha,p}(0,T;B)$ the subset of $L^p(0,T;B)$ function $g$ such that
\begin{align*}
\|g\|_{W^{\alpha,p}(0,T;B)} 
:= \left( \int_0^T \|g(t)\|_B^p dt +
\int_0^T\int_0^T\frac{ \|g(t)-g(s) \|_{B}^p}{|t-s|^{1+\alpha p}}dsdt \right)^{\frac{1}{p}}<\infty.
\end{align*}
\item[(viii)] We denote by $C([0,T]; B)$ the $B$-valued functions that are continuous on $[0,T].$ And for $\alpha>0,$ $C^{\alpha}([0,T];B)$ denotes the set of  $B$-valued $\alpha$-H\"{o}lder  functions $f$ such that 
\begin{align*}
\sup_{t,s\in[0,T],t\neq s}\frac{ \|f(t)-f(s) \|_{B}}{|t-s|^{\alpha}}<\infty.
\end{align*} 
\textcolor{black}{\item[(ix)] if $f$ and $g$ are two quantities, we use $f \lesssim g$ to denote the statement 
that $f \le Cg$ for some constant $C>0$. When this constant $C$ depends on some parameters $a_1,......, a_k,$ 
we use $f \lesssim_{a_1, ..., a_k} g$ to enlighten this dependence on the parameters.}
\end{itemize}
\vspace{3mm}

\vspace{3mm}

Let $\{\beta_{k,l}\}_{k,l\in\mathbb{Z}}$ be a series of independent Brownian motions  on a probability space $(\Omega, \mathcal{F}, \mathbb{P})$.
For $(x, y)\in\mathbb{T}^2, t\in[0, T],$ a $L^2(\mathbb{T}^2)$ cylindrical Wiener process $\xi$ is written by 
$$
\xi(t,x,y)=\sum_{k, l \in\mathbb{Z}}\beta_{k,l}(t)e_{k,l}(x, y).$$
We will see later in Section 5 that the $\varepsilon$-scaled Wiener process is defined as
$$\xi_{\varepsilon}(x,y,t) =  \frac{1}{2L}
\sum_{k,l \in \mathbb{Z}} \beta_{k,l} (\varepsilon^2 t) 
e^{\frac{-\varepsilon i \pi (k,l)}{L} \cdot (x,y)}. 
$$

Here we note the propositions which will be useful later. 

\begin{prop}[Compact embedding 1] \label{CE1} Let $B_0\subset B\subset B_1$ be Banach spaces, $B_0$ and $B_1$ reflexive, with compact embedding of $B_0$ in $B$. Let $p\in(1,\infty)$ and $\alpha\in(0,1)$ be given. Let $X$ be the space 
\begin{equation*}
X=L^p(0,T;B_0)\cap W^{\alpha,p}(0,T; B_1). 
\end{equation*}
Then the embedding of $X$ in $L^p(0,T;B)$ is compact.
\end{prop}
\begin{pr}

See Lemma 2.1 in \cite{FG}.
\end{pr}


\begin{prop}[Compact embedding 2]\label{CE2} If $B_1\subset \tilde{B}$ are two Banach spaces with compact embedding, and the real number $\alpha\in(0,1), p>1$ satisfy 
\begin{align*}
\alpha p>1
\end{align*}
then the space $W^{\alpha,p}(0,T;B_1)$ is compactly embedded into $C([0,T];\tilde{B})$.
\end{prop}
\begin{pr}

See Theorem 2.2 in \cite{FG}.
\end{pr}






\begin{prop}[Gy\"{o}ngy-Krylov criterion] \label{prop:gk}
Let $(X_n)_{n\in\mathbb{N}}$ be a sequence of random variables from a probability space $(\Omega,\mathcal{F},\mathbb{P})$ to a complete separable metric space $(E,d)$. Assume that, for every pair of subsequences $\left(n_1(k),n_2(k)\right),$ with $n_1(k)\ge n_2(k)$ for every $k\in\mathbb{N},$ there is a subsequence $(k(h))_{h\in\mathbb{N}}$ such that the random variables $\left(X_{n_1(k(h))}, X_{n_2(k(h))}\right)$ from $(\Omega,\mathcal{F},\mathbb{P})$ to $(E \times E, d \times d)$ converge in law to a measure $\mu$ on $E \times E$ such that $\mu(D)=\mu(\{(x,y)\in E \times E;x=y\})=1$. Then there exists a random variable $X$ from $(\Omega,\mathcal{F},\mathbb{P})$ to $(E,d)$ such that $X_n$ converges to $X$ in probability.
\end{prop}
\begin{pr}

See Lemma 9 in \cite{FFL}, and \cite{GK}.
\end{pr}

\subsection{Main Theorem}

First of all we set $\varepsilon=1$ and we establish the existence of a solution of  the stochastic Swift-Hohenberg equation:
\begin{equation}\label{SSH}
\begin{cases}
\partial_tu=\mathcal{L}_{1}u-u^3+\dot{\xi}, \,\,\, &{\rm on}\,\, [0,T]\times\mathbb{T}^2,\\
u(0)=u_0,\,\,\,&{\rm on}\,\, \mathbb{T}^2.
\end{cases}
\end{equation}
To find a solution of (\ref{SSH}), we use the decomposition $v=u-Z$ with $Z$ satisfying 
\begin{equation}\label{OS}
\begin{cases}
\partial_tZ=\mathcal{L}_{1}Z+\dot{\xi}, \,\,\, &{\rm on}\,\, [0,T]\times\mathbb{T}^2,\\
Z(0)=0,\,\,\,&{\rm on}\,\, \mathbb{T}^2.
\end{cases}
\end{equation} 
Then $v$ satisfies the following equation formally.
\begin{equation}
\begin{cases}\label{DE}
\partial_tv=\mathcal{L}_{1}v-(v+Z)^3, \,\,\, &{\rm on}\,\, [0,T]\times\mathbb{T}^2,\\
v(0)=u_0,\,\,\,&{\rm on}\,\, \mathbb{T}^2.
\end{cases}
\end{equation}
\textcolor{black}{This equation is a random PDE. We can thus solve the equation $(\ref{DE})$ as a deterministic PDE. }As a result, we can get the solution of $(\ref{SSH}).$

Our main results are as follows.
\begin{thm} \label{ZZ}
Let $T>0$ be fixed. Let $p \ge 1$, $s\in [0,\frac18)$, and $\alpha \in (0, \frac18-s).$ The solution $Z$ 
of (\ref{OS}) has a modification in $C^{\alpha} ([0,T]; {\mathcal{W}}^{s,p}(\mathbb{T}^2))$. Moreover, there exists 
a positive constant $M_{p,T,L}$ such that 
 $$\E\left(\sup_{t \in [0,T]} \|Z(t)\|_{L^p(\mathbb{T}^2)} \right) \le M_{p,T,L}.$$
\end{thm}
\begin{thm}\label{V}
Let $u_0\in L^2(\mathbb{T}^2),$ and $\alpha>0$. There exists a \textcolor{black}{unique} stochastic process $v$ on $(\Omega,\mathcal{F},\mathbb{P})$ satisfying 
(\ref{SSHE}) in the following weak sense, i.e.,  
for $w\in \mathcal{W}^{1+\alpha,2}$, $t\in[0,T]$,
\begin{align}
(v_t,w)=(u_0,w)+\int_0^t(v_s,\mathcal{L}_{1}w)ds-\int_0^t((v_s+Z_s)^3,w)ds,
\end{align}
and $v$ takes values in $L^{\infty}(0,T;L^2(\mathbb{T}^2))\cap L^3(0,T;L^4)\cap C([0,T];\mathcal{W}^{-(1+\alpha),2}) \cap L^2(0,T; \mathcal{W}^{1-\alpha,2})$ almost surely.
\end{thm}

Theorem $\ref{ZZ}$ is proved by using the Kolmogorov test, where convergence properties of the Gamma function are helpful in the computation. 
Theorem $\ref{V}$ is proved by using a Galerkin approximation as in \cite{BS}. 
Note that $L^p$ energy estimates seem not available, thus we do not use the fixed point argument. 
First we consider a finite dimensional nonlinear equation. We get an energy estimate and properties of $Z$ allow us to obtain a probabilistically uniform estimate with respect to the dimension. The Prohorov Theorem and the Skorohod Theorem imply the existence of a limit 
taking a subsequence on another probability space. Moreover, the Gy\"ongy-Krylov criterion can make the convergence on another probability space into the convergence on the original space regarding $X_n$ as the subsequence converging to some probability measure weakly. This convergence constructs a solution of the infinite dimensional system.


\section{Regularity of the solution of the stochastic linear equation}

\subsection{Regularity of $Z$}  In this section, we investigate the regularity of $Z$ by using the Kolmogorov test.
We may write $Z$ as a mild solution.
\begin{equation}\label{Z}
Z(t)=\int_0^te^{(t-s)\mathcal{L}_{1}}d\xi(s)
=\sum_{k,l}\int_0^te^{(t-s)\lambda_{k,l,1}}e_{k,l} (x,y) d\beta_{k,l}(s).
\end{equation}

\begin{prop} \label{prop:Kolmo}
Let $T>0$ be fixed. Let $p \ge 1$, and $0 \le \theta<\frac{1}{8}$. Then $Z$ has a modification which is 
$\alpha$-$H\ddot{o}lder$ continuous on $[0,T]$ with values in ${\mathcal{W}}^{\theta,p}(\mathbb{T}^2)$ 
for $\alpha \in (0, \frac18-\theta)$.
\end{prop}
\begin{pr}

Let $(x,y)\in\mathbb{T}^2$ and $\theta \ge 0$. For $t>s$, $t,s \in [0,T]$, we first calculate
\begin{align*}
&\mathbb{E} \left( |(1-\mathcal{L}_0)^{\frac{\theta}{2}} (Z(t) -Z(s)) |^2 \right) \\
=&\mathbb{E} \sum_{k,l \in \Z} ( (1- {\textcolor{black}{ \left(\frac {\pi}{L} k \right)}}^2)^{2}
+ {\textcolor{black}{\left(\frac{\pi}{L} l \right)}}^{2}+1)^{\theta} (Z(t,x,y)-Z(s,x,y),e_{k,l})^2 \\
=&\sum_{k,l}\mathbb{E}\left|\int_0^s((1-\textcolor{black}{\left(\frac {\pi}{L} k\right)}^2)^2
+ {\textcolor{black}{\left(\frac{\pi}{L} l \right)}}^2+1)^{\theta}(e^{(t-u)\lambda_{k,l,1}}-e^{(s-u)\lambda_{k,l,1}}) 
d\beta_{k,l}(u)\right|^2\\
+&\sum_{k,l}\mathbb{E}\left|\int_s^t((1-\textcolor{black}{\left(\frac {\pi}{L} k\right) }^2)^2
+ {\textcolor{black}{\left(\frac{\pi}{L} l\right)}}^2+1)^{\theta}e^{(t-u)\lambda_{k,l,1}}  d\beta_{k,l}(u)\right|^2\\
\le&\sum_{k,l}\int_0^s((1-\textcolor{black}{\left(\frac {\pi}{L} k\right)}^2)^2
+ {\textcolor{black}{\left(\frac{\pi}{L} l\right)}}^2+1)^{2\theta}(e^{(t-u)\lambda_{k,l,1}}-e^{(s-u)\lambda_{k,l,1}})^2du\\
+&\sum_{k,l}\int_s^t((1-\textcolor{black}{\left(\frac {\pi}{L} k\right)}^2)^2
+ {\textcolor{black}{\left(\frac{\pi}{L} l\right)}}^2+1)^{2\theta}e^{2(t-u)\lambda_{k,l,1}}du\\
=&I_1+I_2,
\end{align*}
where we have used the It\^{o} isometry of the stochastic integral.
First, we estimate $I_2$ dividing into $I_{2,\ge 0}$ and $I_{2, <0}$ with   
\begin{eqnarray*}
I_{2,\ge 0} &=& \sum_{ \substack{l\in \Z,\\ k: 1- k^2 \ge 0}} \int_s^t((1-\textcolor{black}{\left(\frac{\pi}{L} k \right)^2)}^2+\textcolor{black}{\left(\frac{\pi}{L}l\right)^2}+1)^{2\theta} 
e^{2(t-u)\lambda_{k,l,1}}du,\\
I_{2, <0} &=& \sum_{\substack{ l\in \Z,\\ k: 1-k^2 < 0}} \int_s^t((1-\textcolor{black}{\left(\frac{\pi}{L} k\right)^2})^2+ \textcolor{black}{\left(\frac{\pi}{L} l \right)^2}+1)^{2\theta} 
e^{2(t-u)\lambda_{k,l,1}}du.
\end{eqnarray*}
Recall $\lambda_{k,l,1}=-\{(1-\textcolor{black}{\left(\frac {\pi}{L} k\right)}^2)^2 + {\textcolor{black}{\left(\frac{\pi}{L} l\right)}}^2 -1\}$. Note that if $k\in \Z$ 
satisfies $1-k^2 \ge 0$, then $k=0,\pm1$ thus $(1-\textcolor{black}{\left(\frac {\pi}{L} k\right)}^2)^2 \le 1$. Therefore, 
\begin{eqnarray*} 
I_{2, \ge 0} &\le & \sum_{\substack{ l\in \Z,\\ k: 1-k^2 \ge  0}} 
\int_s^t ( {\textcolor{black}{\left(\frac{\pi}{L} l\right)}}^2+2)^{2\theta}e^{2(t-u)\lambda_{k,l,1}}du \\
&\le & \textcolor{black}{\frac{4L}{\pi}} \int_s^t \int_0^{+\infty} \int_0^1 (2+y^2)^{2 \theta} e^{ -(t-u)\{(1-\textcolor{black}{(\frac{\pi}{L})^2} x^2)^2+y^2-1\}} dx dy du \\
& \le & \textcolor{black}{\frac{4L}{\pi}} \int_s^t \int_0^{+\infty} (2+y^2)^{2 \theta} e^{ -(t-u)(y^2-1)} dy du \\
&=& \textcolor{black}{\frac{4L}{\pi}} \int_s^t e^{(t-u)} \int_0^{+\infty}   (2+y^2)^{2 \theta} e^{ -(t-u)y^2} dy du
\end{eqnarray*}
Here, the use of the change of variable $(t-u) y^2=z \ge 0$ allows us to estimate
\begin{eqnarray*}
\int_0^{+\infty}   (2+y^2)^{2 \theta} e^{ -(t-u)y^2} dy 
&\lesssim_{\theta}& \int_0^{+\infty} e^{ -(t-u)y^2} dy + \int_0^{+\infty}   y^{4 \theta} e^{-(t-u)y^2} dy \\
&=& C(\theta) (t-u)^{-\frac12} +C'(\theta) (t-u)^{-2\theta-\frac12} \Gamma \left(2\theta+\frac12\right)
\end{eqnarray*} 
if $2\theta+\frac12>0$, where $\Gamma(\cdot)$ is the Gamma function, and 
$$\Gamma \left(2\theta+\frac12\right)=\int_0^{+\infty} z^{2\theta-\frac12} e^{-z} dz.$$
Hence, 
\begin{eqnarray*}
I_{2, \ge 0} &\lesssim_{\theta, \textcolor{black}{L}} & \int_s^t e^{(t-u)} \left( (t-u)^{-\frac12} +(t-u)^{-2\theta-\frac12} \right) du 
\lesssim_{\theta, T, \textcolor{black}{L}} (t-s)^{-2\theta +\frac12},
\end{eqnarray*}
if $\theta \in [0, \frac14)$. Next we estimate $I_{2,<0}$. 
\textcolor{black}{ 
The condition is $k^2>1$, but thanks to the symmetry, 
we first focus on the integral
\begin{eqnarray} \label{largerthan1}
\int_{\frac{\pi}{L}}^{\infty} (1-x^2)^{4\theta} e^{-2(t-u)(1-x^2)^2} dx.
\end{eqnarray}
}
Let $x^2-1=z$, and we get\footnote{\textcolor{black}{considering two cases: $\frac{\pi}{L}<1$ and $x>1$, or $\frac{\pi}{L} \ge 1$ or $x>1$. The former case is impossible in (\ref{largerthan1}). }}
$$
(\ref{largerthan1}) \le \int_1^{\infty} (1-x^2)^{4\theta} e^{-2(t-u)(1-x^2)^2} dx
\lesssim_{\theta} \int_0^{\infty }z^{4\theta-\frac{1}{2}}e^{-2(t-u)z^2}\frac{1}{2(z+1)^{\frac{1}{2}}}dz.
$$
Moreover we change the variable $2(t-u) z^2 = w\ge 0$ which leads to the above RHS:
\begin{eqnarray*}
&\lesssim_{\theta}& \int_0^{\infty} \frac{w^{\frac{1}{2}(4\theta-\frac{1}{2})}}{(t-u)^{\frac{1}{2}(4\theta-\frac{1}{2})}}\frac{e^{-w}}{w^{\frac{1}{2}}(t-u)^{\frac{1}{2}}}dw\\
&\lesssim_{\theta}& (t-u)^{-\frac{1}{4}-2\theta},
\end{eqnarray*}
where we have again used the convergence of the Gamma function if  $2\theta+\dfrac14 >0.$ 
Therefore, by symmetry, 
\begin{eqnarray*}
I_{2, <0} &=& \sum_{\substack{ l\in \Z,\\ k: 1-k^2 < 0}} \int_s^t ((1- {\textcolor{black}{\left(\frac{\pi}{L} k\right)}}^2)^2+ {\textcolor{black}{\left(\frac{\pi}{L} l\right)}}^2+1)^{2\theta} 
e^{2(t-u)\lambda_{k,l,1}}du \\
& \lesssim_{\textcolor{black}{L}} & \int_s^t \int_0^{+\infty} \int_{\textcolor{black}{\frac{\pi}{L}}}^{+\infty} 
 ((1-x^2)^2+y^2+1)^{2\theta} e^{-2(t-u)\{ (1-x^2)^2+y^2-1\}} dxdydu \\
 &\lesssim_{\textcolor{black}{\theta, L}}& \int_s^t \int_0^{+\infty} \int_{\textcolor{black}{\frac{\pi}{L}}}^{+\infty} 
 (1-x^2)^{4\theta} e^{-2(t-u)\{ (1-x^2)^2+y^2-1\}} dxdydu \\
 && \hspace{3mm}+
 \int_s^t \int_0^{+\infty} \int_{\textcolor{black}{ \frac{\pi}{L} }}^{+\infty} 
 (y^{2\theta}+1) e^{-2(t-u)\{ (1-x^2)^2+y^2-1\}} dxdydu \\
&\lesssim_{\textcolor{black}{\theta,L}}&  \int_s^t  e^{2(t-u)} (t-u)^{-\frac14 -2\theta} \int_0^{+\infty}  
 e^{-2(t-u)y^2} dydu \\
 && \hspace{3mm}+ 
 \int_s^t  e^{2(t-u)}  \int_0^{+\infty} \int_1^{+\infty} 
 (y^{2\theta}+1) e^{-2(t-u)\{ (1-x^2)^2+y^2\}} dxdydu \\
&\lesssim_{\theta,T,\textcolor{black}{L}}& \int_s^t (t-u)^{-\frac34 -2\theta} du 
+\int_s^t (t-u)^{-\frac14} \{(t-u)^{-\frac12} +(t-u)^{-\frac12-\theta} \} du
\lesssim_{T,\theta, \textcolor{black}{L} } (t-s)^{\frac14 -2\theta}
\end{eqnarray*}
if $\theta <\dfrac18$ 
which implies,
$I_2\lesssim_{\theta,T, \textcolor{black}{L} }(t-s)^{\frac{1}{4}-2\theta}$ if $\theta<\dfrac18$. Now we estimate $I_1$. 
For $\gamma\in(0,1),$ a similar calculation as above yields, using the $\gamma$-H\"{o}lder regularity of the exp function,   
\begin{eqnarray*}
&&I_1\lesssim \int_0^s\int_{\mathbb{R}^2}(t-s)^{2\gamma}((1- {\textcolor{black}{\left(\frac{\pi}{L}x\right)}}^2)^2+ {\textcolor{black}{\left(\frac{\pi}{L} y\right)}}^2+1)^{2\gamma+2\theta}
e^{-(s-u)\{ (1- {\textcolor{black}{\left(\frac{\pi}{L} x\right)}}^2)^2+ {\textcolor{black}{\left(\frac{\pi}{L} y\right)}}^2-1\} }dxdydu\\
&&\lesssim_{\textcolor{black}{\gamma, \theta}} (t-s)^{2\gamma} \int_0^s\int_{\mathbb{R}^2} 
(t-s)^{2\gamma}((1- {\textcolor{black}{\left(\frac{\pi}{L} x\right)}}^2)^{4(\gamma+\theta)}
+ {\textcolor{black}{\left(\frac{\pi}{L} y\right)}}^{4(\gamma+\theta)}+1)
e^{-(s-u) \{ (1- {\textcolor{black}{\left(\frac{\pi}{L} x\right)}}^2)^2+ {\textcolor{black}{\left(\frac{\pi}{L} y\right)}}^2-1\}}dxdydu\\
&&\lesssim_{\textcolor{black}{\gamma,\theta,T,L}} (t-s)^{2\gamma}\int_0^s(s-u)^{-\frac{3}{4}-2\gamma-2\theta}du.
\end{eqnarray*}
The right hand side is finite if $\frac{1}{8}>\theta+\gamma$. Hence, for $m\in\mathbb{N}$, we obtain
\begin{align*}
\mathbb{E} \left( |(1-\mathcal{L}_0)^{\frac{\theta}{2}} (Z(t) -Z(s)) |^{2m} \right)
\le C(m,\gamma,\theta,T)|t-s|^{\min(2\gamma,\frac{1}{4}-2\theta)\times 2m}. 
\end{align*}
Therefore, for $m\in\mathbb{N}$ and $1\le p\le 2m,$ by the Minkowski inequality, 
\begin{eqnarray*}
\E\left( \|(1-\mathcal{L}_0)^{\frac{\theta}{2}} (Z(t) -Z(s)) \|_{L^p(\mathbb{T}^2)}^{2m} \right)^{\frac{1}{2m}} 
&\le& C(m,\gamma,\theta,T)|t-s|^{\gamma}.\\
\end{eqnarray*}
Set $\theta=0$. We conclude by the Kolmogorov test 
that $Z$ has a modification in $C^{\alpha}([0,T];\mathcal{W}^{\theta,p}(\mathbb{T}^2))$ 
for any $\alpha <\frac18$ and $p \ge 1$. In particular, 
 $$\E\left(\sup_{t \in [0,T]} \|Z(t)\|_{L^p(\mathbb{T}^2)} \right) \le M_{p,T,L},$$
 for some $M_{p,T,L}>0.$ 
More generally, if $\theta \in (0,\frac18)$, $Z$ has a modification in $C^{\alpha}([0,T];\mathcal{W}^{\theta,p}(\mathbb{T}^2))$ for 
$p \ge 1$, and $\alpha < \frac18-\theta$.    
\end{pr}

\section{Existence of the solution}
In this section, we construct a solution of $(\ref{DE})$ using a compactness argument.

\subsection{Approximation}
We want to construct the solution of ($\ref{DE}$) by a Faedo Galerkin approximation.
For $f\in L^2(\mathbb{T}^2)$, and $n\in\mathbb{N},$ 
\begin{align*}
\Pi_n:L^2(\mathbb{T}^2)\rightarrow \Pi_nL^2(\mathbb{T}^2)
\end{align*}
is defined by
\begin{equation}
\Pi_nf=\sum_{|k|+|l|\le n}(e_{k,l},f)e_{k,l},
\end{equation}
where $\Pi_nL^2(\mathbb{T}^2)$ denotes the subspace of $L^2(\mathbb{T}^2)$ such that the $\Pi_n f$ can be represented by linear combinations of $e_{k,l}$ with $|k|+|l|\le n$. 
Obviously, $\Pi_nf\rightarrow f$ in $L^2(\mathbb{T}^2)$.
Our goal is to find the solution of $(\ref{DE})$.  It is defined by 
\begin{equation}
(v(t)-u_0,w)=\int_0^t (v(\sigma),\mathcal{L}_{1}w)d\sigma-\int_0^t ((v(\sigma)+Z(\sigma))^3,w)d\sigma,\,\,\,\,
\mathbb{P}-a.s.
\end{equation}
for  $w\in \mathcal{W}^{1+\alpha,2}$. 
To do so we use the finite dimensional solution $v_n$ which solves 
\begin{equation}{\label{FDE}}
\begin{cases}
\partial_tv_n(t)&=\mathcal{L}_{1} v_n(t)-\Pi_n(\Pi_nv_n(t)+\Pi_nZ(t))^3,\\
v_n(0)&=\Pi_nu_0.
\end{cases}
\end{equation}
This solution \textcolor{black}{is smooth enough to obtain} the following energy estimate.
Multiplying the equation $(\ref{FDE})$ by $v_n$, and using periodic boundary condition, we get:
\begin{align*}
\frac{1}{2}\frac{d}{dt}\|v_n\|_{L^2}^2
&=-2\int_{\mathbb{T}^2}v_n\partial_x^2v_ndxdy\\
&-\int_{\mathbb{T}^2}v_n \partial_x^4 v_ndxdy+\int_{\mathbb{T}^2} v_n \partial_y^2vdxdy-\int_{\mathbb{T}^2}v_n(v_n+\Pi_nZ)^3dxdy.
\end{align*}
The second term is estimated as :
\begin{align*}
\int_{\mathbb{T}^2}v_n\partial_x^2v_ndxdy\le \int_{\mathbb{T}^2}\frac{v_n^2}{2\delta^2}dxdy+\int_{\mathbb{T}^2}\frac{\delta^2(\partial_x^2v_n)^2}{2}dxdy,\,\,\,(\delta>0), 
\end{align*}
and the fourth term may be written as, by integration by parts, 
\begin{align*}
\int_{\mathbb{T}^2}v_n\partial_x^4v_ndxdy
=-\int_{\mathbb{T}^2}(\partial_xv_n)(\partial_x^3v_n)dxdy
=\int_{\mathbb{T}^2}(\partial_x^2v_n)^2dxdy
=\|\partial_x^2v_n\|_{L^2}^2,
\end{align*}
and 
\begin{align*}
\int_{\mathbb{T}^2}v_n(v_n+\Pi_nZ)^3dxdy \le -\frac{1}{2} \|v_n \|_{L^4}^4+C\|\Pi_nZ \|_{L^4}^4,
\end{align*}
where we have used Young's inequality.
Then for any $\delta\in(0,1)$,
\begin{align*}
-2\int_{\mathbb{T}^2}v_n\partial_x^2v_ndxdy-\int_{\mathbb{T}^2}v_n \partial_x^4v_ndxdy+\int_{\mathbb{T}^2}v_n\partial_y^2v_ndxdy-\int_{\mathbb{T}^2}v_n(v_n+\Pi_nZ)^3dxdy\\
<\frac{1}{\delta^2} \|v_n \|_{L^2}^2+(\delta^2-1) \|\partial_x^2v_n \|_{L^2}^2
-\|\partial_yv_n\|_{L^2}^2-\frac{1}{2} \|v_n \|_{L^4}^4+C\|\Pi_nZ \|_{L^4}^4.
\end{align*}
Therefore, we have,  
\begin{equation*}
\frac{d}{dt} \|v_n \|_{L^2}^2+ \|\partial_x^2v_n \|^2_{L^2}+ \|\partial_yv_n \|_{L^2}^2
+\|v_n \|_{L^4}^4
\le C(\delta) \left( \|v_n \|_{L^2}^2+ \|\Pi_nZ \|_{L^4}^4\right),
\end{equation*}
which reads in other word, by the definition of the space $\mathcal{W}^{1,2}$, taking for example $\delta=\dfrac12$,
\begin{equation} \label{EE}
\frac{d}{dt} \| v_n \|_{L^2}^2+\|v_n \|_{\mathcal{W}^{1,2}}^2+\|v_n \|_{L^4}^4
<C\left( \|v_n \|_{L^2}^2+ \|\Pi_nZ \|_{L^4}^4\right).
\end{equation}
Thus, by integrating on $[0,t]$ with $t \le T$,  
\begin{eqnarray}  \nonumber
&& \|v_n\|_{L^2}^2 + \int_0^t \|v_n(s)\|_{\mathcal{W}^{1,2}}^2 ds + \int_0^t  \|v_n(s) \|_{L^4}^4 ds \\ \label{eq:integral}
&\le& C\int_0^t (\|v_n(s)\|_{L^2}^2 + \|\Pi_n Z\|_{L^4}^4) ds +\| \Pi_n u_0\|_{L^2}^2.  
\end{eqnarray}      
Therefore, by the Gronwall inequality, we have  
\begin{equation}\label{L2 of v}
\|v_n \|_{L^2}^2(t) \le C e^{CT} \left(\|\Pi_nu_0 \|_{L^2}^2
+\int_0^T \|\Pi_nZ \|_{L^4}^4(s)ds\right).
\end{equation}
On the other hand, a similar proof as in Proposition \ref{CE2} 
infers that $\mathbb{E} \|\Pi_nZ \|_{L^4}^4\le C_T$ where $C_T$ is independent of $n$.
This implies that
\begin{align} \nonumber
\mathbb{E} \|v_n \|_{L^2}^2 &\le C(\varepsilon) \|\Pi_n u_0 \|_{L^2}^2 e^{CT}
+C\int_0^T\mathbb{E} { \|\Pi_nZ \|_{L^4}^4}ds\\ \label{EE2n}
&\le C(T) \|\Pi_nu_0 \|_{L^2}^2+C(T).
\end{align}
Thus, by (\ref{eq:integral}) 
\begin{equation} \label{EEE2n}
\mathbb{E}\left(\int_0^T \|v_n \|_{\mathcal{W}^{1,2}}^2ds+\int_0^T \|v_n \|_{L^4}^4ds\right)
\le C(T, \|\Pi_nu_0 \|_{L^2}).
\end{equation}
The fact $\|\Pi_nu_0 \|_{L^2} \le \|u_0 \|_{L^2}$ implies that $\{v_n\}_n$ 
is bounded (independently of $n$) in  
\begin{align*}
L^2(\Omega; L^2(0,T;\mathcal{W}^{1,2}))\cap L^4(\Omega;L^4(0,T; L^4(\mathbb{T}^2))). 
\end{align*}
Now for any $w\in \mathcal{W}^{1,2}$,  noting that $\mathcal{L}_1 = 2-(1- \mathcal{L}_0)$, 
\begin{align*}
\left|\left(\frac{d}{dt}v_n(t),w\right)\right|&\le|(v_n(t),\mathcal{L}_{1}w)|+|((v_n+\Pi_nZ)^3,\Pi_nw)|\\
&\le \|v_n(t) \|_{\mathcal{W}^{1,2}} \|w\|_{\mathcal{W}^{1,2}}+\|(v_n+\Pi_nZ)^3\|_{L^{\frac{4}{3}}} \|\Pi_nw \|_{L^4}\\
&\le C\left( \|v_n(t) \|_{\mathcal{W}^{1,2}} \|w\|_{\mathcal{W}^{1,2}}+(\|v_n\|_{L^4}^3+\|\Pi_nZ\|_{L^4}^3)\|w\|_{\mathcal{W}^{1,2}}\right).
\end{align*}
Therefore, 
\begin{align*}
\left|\left|\frac{d}{dt}v_n(t)\right|\right|_{\mathcal{W}^{-1,2}}
&\le C\left(\|v_n(t) \|_{\mathcal{W}^{1,2}}+ \|v_n \|_{L^4}^3+ \|\Pi_nZ \|_{L^4}^3\right).\\
\end{align*}
The uniform estimates of  (\ref{EEE2n}) and $ \|\Pi_nZ \|_{L^4}^4$ in $n$ imply
\begin{equation}
\mathbb{E}\int_0^T\left|\left|\frac{d}{dt}v_n(t)\right|\right|_{\mathcal{W}^{-1,2}}^{\frac{4}{3}}dt
\le C(T, \|u_0 \|_{L^2}),
\end{equation}
which concludes that  $\{v_n\}$ is bounded 
in $L^{\frac{4}{3}}(\Omega;W^{1,\frac{4}{3}}(0,T; \mathcal{W}^{-1,2}))$. 
Let us now consider $m>0$. We multiply both sides of  (\ref{EE}) by $\|v_n\|_{L^2}^{2m-2}$,
\begin{equation} \label{Energym}
\frac{1}{m}\frac{d}{dt} \|v_n \|_{L^2}^{2m} +\|v_n \|_{\mathcal{W}^{1,2}}^2 \|v_n \|_{L^2}^{2m-2}
\le C(\|v_n \|_{L^2}^{2m}+ \|\Pi_nZ \|_{L^4}^4 \|v_n \|_{L^2}^{2m-2}).
\end{equation}
Applying the interpolation inequality:
$$
\| v_n \|_{\mathcal{W}^{\frac23,2}} 
\le \|v_n\|_{L^2}^{1/3}  \|v_n \|_{\mathcal{W}^{1,2}}^{2/3}
$$
and choosing $m=3$ in (\ref{Energym}), we want to show that $\{v_n\}$ 
is bounded in $L^3(\Omega; L^3(0,T; \mathcal{W}^{\frac23,2}))$. 

We integrate (\ref{Energym}) in time, and we get
\begin{eqnarray*}
&&\frac13 \|v_n\|_{L^2}^6 +3 \int_0^T \|v_n\|_{\mathcal{W}^{\frac23,2}}^{12} \\
&&  \le \frac13 \|v_0\|_{L^2}^6 + C\int_0^T \left( \|v_n\|_{L^2}^6 +\|\Pi_n Z\|_{L^4}^4 \|v_n\|_{L^2}^4 \right) \\ 
&& \le \frac13 \|v_0\|_{L^2}^6 +C\int_0^T \|v_n\|_{L^2}^6 + C\int_0^T \|\Pi_n Z\|_{L^4}^{12}. 
\end{eqnarray*}
We take the expectation and the use of the Gronwall inequality \textcolor{black}{which} implies that 
\begin{eqnarray*}
\E(\|v_n\|_{L^2}^6) &\lesssim& \|u_0\|_{L^2}^6 + \E \left(\int_0^T \|\Pi_n Z\|_{L^4}^{12} \right) \\
&\lesssim& \|u_0\|_{L^2}^6 + \left(\int_{\mathbb{T}^2} \E(|\Pi_n Z|^{12})^4 dxdy \right)^3.
\end{eqnarray*}
Here we have used the Minkowski inequality, and \textcolor{black}{that} the integrand of the second term in the RHS may be bounded 
by $C_T$ as in the proof of Proposition \ref{prop:Kolmo}. 
Therefore, we obtain 
\begin{align*}
\mathbb{E}\left(\int_0^T \|v_n \|_{\mathcal{W}^{\frac{2}{3},2}}^3 ds\right)\le C(T, \|u_0\|_{L^2}).
\end{align*}
We remark that embedding $L^3(0,T; \mathcal{W}^{\frac23,2} )\cap W^{1,\frac{4}{3}} (0,T; \mathcal{W}^{-1,2})
\subset L^{3}(0,T; L^4)$ is compact , and $L^2 (0,T; \mathcal{W}^{1,2} )\cap W^{1,\frac{4}{3}} (0,T; \mathcal{W}^{-1,2}) 
\subset L^{2}(0,T; \mathcal{W}^{1-\alpha,2})$ is compact for any $\alpha>0$ 
by Proposition $\ref{CE1}$. 
On the other hand the embedding $W^{1,\frac{4}{3}}(0,T; \mathcal{W}^{-1,2} )\subset C([0,T]; \mathcal{W}^{-(1+\alpha),2})$ is compact for any $\alpha>0$ by Proposition $\ref{CE2}$.
Meanwhile, as in the proof of  Proposition \ref{prop:Kolmo},  
$\{\Pi_nZ\}$ is bounded in 
$C^{\alpha_0}([0,T]; \mathcal{W}^{s,p})$ for $0 \le s<\frac{1}{8}$, $\alpha_0< \frac{1}{8} -s$, $p \ge 1$.  
The embedding $C([0,T]; \mathcal{W}^{\frac{1}{16},4})\cap C^{\alpha_0}([0,T]; \mathcal{W}^{-\delta,4}(\mathbb{T}^2)) 
\subset C([0,T];L^4)$ is thus compact for any $\delta>0$. 
We thus conclude that the sequence $\{(v_n, \Pi_nZ)\}_{n}$ is tight in $L^3(0,T; L^4)\cap C([0,T]; \mathcal{W}^{-(1+\alpha),2}) \cap  L^{2}(0,T; \mathcal{W}^{1-\alpha,2})  \times C([0,T];L^4)$ for any $\alpha>0$.

\subsection{Existence of the solution}

By Prohorov Theorem, 
the tightness of $(v_n,\Pi_n Z)$ implies that the existence of subsequence $(v_{n(k)},\Pi_{n(k)} Z)$ and 
some probability measure $\mu$ such that 
\begin{align*}
(v_{n(k)},\Pi_{n(k)} Z)\rightarrow \mu\,\,\,{\rm weakly~ as}~  k\rightarrow\infty. 
\end{align*}
Moreover by Skorohod Theorem, there exists a probability space $(\Omega',\mathcal{F}',\mathbb{P}')$ and random variables $\{X^{n(k)},Z^{n(k)}\}_k$, and $(X,Z)$ such that 
\begin{align*}
& \mathrm{Law}(v_{n(k)},\Pi_{n(k)}Z)=\mathrm{Law}(X^{n(k)},Z^{n(k)})\,\,\,\,\,\mbox{for} ~ k \ge 1, \quad  \mathrm{Law}(X,Z)=\mu,\\
& \lim_{k\rightarrow\infty}(X^{n(k)},Z^{n(k)})=(X,Z), \quad {\rm in}\\ 
& L^3(0,T;L^4)\cap C([0,T]; \mathcal{W}^{-(1+\alpha),2} \cap L^2(0,T; \mathcal{W}^{1-\alpha,2}) \times C ([0,T], L^4),\,\,\mathbb{P}'-a.s.
\end{align*}
The equivalence of probability laws leads that for $w\in \mathcal{W}^{1+\alpha,2}$ and $t\in[0,T]$,
\begin{equation} \label{weak}
\left(X^{n(k)}(t)-u_{n(k)}(0), w\right)=\int_0^t\left(X^{n(k)}(\sigma),\mathcal{L}_{1}w\right)d\sigma-\int_0^t\left(\Pi_{n(k)}(X^{n(k)}+Z^{n(k)})^3(\sigma),w\right)d\sigma
\end{equation}
$\mathbb{P}'$-almost surely.
We prove that $X$ is a solution of $(\ref{DE})$ on $(\Omega', \mathcal{F}', \prob' )$.
For $t\in[0,T],$
\begin{align*}
\left|(X^{n(k)}(t)-X(t),w)\right|\le \|X^{n(k)}-X \|_{C([0,T];\mathcal{W}^{-(1+\alpha),2})} \|w \|_{\mathcal{W}^{1+\alpha,2}}\rightarrow0
\end{align*}
as $k\rightarrow\infty$. The LHS of (\ref{weak}) converges to $\left(X(t)-u_0,w\right)$.
Next we observe the convergence of the second term on the RHS: 
\begin{align*}
\left|\int_0^t\left(X^{n(k)}(\sigma)-X(\sigma),\mathcal{L}_{1} w\right)d\sigma\right|&\le
\int_0^t  \left( X^{n(k)}(\sigma)-X(\sigma), (2-(1-\mathcal{L}_0)) w \right) d\sigma\\
&\le\int_0^T \|X^{n(k)}(\sigma)-X(\sigma) \|_{\mathcal{W}^{-1-\alpha,2}} \|2w\|_{\mathcal{W}^{1+\alpha,2}}d\sigma\\
&+ \int_0^T \|X^{n(k)}(\sigma)-X(\sigma) \|_{\mathcal{W}^{1-\alpha,2}} \|w\|_{\mathcal{W}^{1+\alpha,2}} d\sigma\\
&\rightarrow 0,
\end{align*}
as $k\rightarrow\infty$.
The convergence of the third term on the RHS can be shown as follows;
\begin{align*}
&\int_0^t\left|\left(\Pi_{n(k)}\left(X^{n(k)}(\sigma)+Z^{n(k)}(\sigma)\right)^3,w\right)-\left(\left(X(\sigma)+Z(\sigma)\right)^3,w\right)\right|d\sigma\\
\le&\int_0^t\left|\left(\left(X^{n(k)}(\sigma)+Z^{n(k)}(\sigma)\right)^3,\Pi_{n(k)}w-w\right)\right|d\sigma\\
+&\int_0^t\left|\left(\left(X^{n(k)}(\sigma)+Z^{n(k)}(\sigma)\right)^3-\left(X(\sigma)+Z(\sigma)\right)^3,w\right)\right| d\sigma\\
=&I+J.
\end{align*}
First, we estimate $I$.
\begin{align*}
I&\le C\int_0^T\|\left(X^{n(k)}+Z^{n(k)}\right)^3 \|_{L^{\frac{4}{3}}} \|\Pi_{n(k)}w-w\|_{L^4} d\sigma \\
&\le C\int_0^T\left(\|X^{n(k)}(\sigma)\|_{L^4}^3
+\|Z^{n(k)}(\sigma)\|_{L^4}^3\right)||\Pi_{n(k)}w-w||_{L^4} d\sigma \\
&\le C\int_0^T\left(\|X^{n(k)}(\sigma)\|_{L^4}^3+\|Z^{n(k)}(\sigma)\|_4^3 \right) 
\|\Pi_{n(k)}w-w \|_{\mathcal{W}^{1+\alpha,2}} d\sigma\\
&\rightarrow 0
\end{align*}
as $k\rightarrow\infty,$ for any $\alpha>0$, where we have used the Sobolev embedding in the third inequality.  
Next we will see the convergence of $J$.
\begin{align*}
J&=\int_0^t\left|\left(\left(X^{n(k)}(\sigma)+Z^{n(k)}(\sigma)\right)^3-\left(X(\sigma)+Z(\sigma)\right)^3,w\right)\right|d\sigma\\
&\le C\int_0^t\left|\left(\left(X^{n(k)}(\sigma)-X(\sigma)\right)\left(\left|X^{n(k)}(\sigma)\right|^2+\left|X(\sigma)\right|^2
+\left|Z^{n(k)}(\sigma)\right|^2+\left|Z(\sigma)\right|^2\right),|w|\right)\right|d\sigma\\
&+C\int_0^t\left|\left(\left(Z^{n(k)}(\sigma)-Z(\sigma)\right)\left(\left|X^{n(k)}(\sigma)\right|^2+\left|X(\sigma)\right|^2+\left|Z^{n(k)}(\sigma)\right|^2
+\left|Z(\sigma)\right|^2\right),|w|\right)\right|d\sigma.\\
&=J_1+J_2.\\
\end{align*}
Considering the first term $J_1$ we have;
\begin{align*}
J_1&\le \int_0^T\left|\left|X^{n(k)}(\sigma)-X(\sigma)\right|\right|_{L^4}\left(\left|\left|X^{n(k)}(\sigma)\right|\right|_{L^4}^2
+\left|\left|X(\sigma)\right|\right|_{L^4}^2+\left|\left|Z^{n(k)}(\sigma)\right|\right|_{L^4}^2+\left|\left|Z(\sigma)\right|\right|_{L^4}^2\right)d\sigma \\
& \times \|w\|_{L^4}\\
&\le \textcolor{black}{C} \left(\int_0^T\left|\left|X^{n(k)}(\sigma)-X(\sigma)\right|\right|_{L^4}^3d\sigma\right)^{\frac{1}{3}} \\
&\times \left(\int_0^T\left(\left|\left|X^{n(k)}(\sigma)\right|\right|_{L^4}^3+\left|\left|X(\sigma)\right|\right|_{L^4}^3+\left|\left|Z^{n(k)}(\sigma)\right|\right|_{L^4}^3+\left|\left|Z(\sigma)\right|\right|_{L^4}^3\right)d\sigma\right)^{\frac{2}{3}} \|w\|_{L^4}\\
\end{align*}
The second integral is finite and the convergence $X^{n(k)}\rightarrow X$ {\rm in} $L^3(0,T; L^4)$ gives 
the convergence of $J_1$. Furthermore, similarly as above, we estimate   
\begin{align*}
J_2&\le \int_0^T\left|\left| Z^{n(k)}(\sigma)-Z(\sigma)\right|\right|_{L^4}\left(\left|\left|X^{n(k)}(\sigma)\right|\right|_{L^4}^2
+\left|\left|X(\sigma)\right|\right|_{L^4}^2+\left|\left|Z^{n(k)}(\sigma)\right|\right|_{L^4}^2+\left|\left|Z(\sigma)\right|\right|_{L^4}^2\right)d\sigma \\
& \times \|w\|_{L^4}\\
&\le \textcolor{black}{C}\left(\int_0^T\left|\left| Z^{n(k)}(\sigma)-Z(\sigma)\right|\right|_{L^4}^3d\sigma\right)^{\frac{1}{3}} \\
&\left(\int_0^T\left(\left|\left|X^{n(k)}(\sigma)\right|\right|_{L^4}^3+\left|\left|X(\sigma)\right|\right|_{L^4}^3+\left|\left|Z^{n(k)}(\sigma)\right|\right|_{L^4}^3+\left|\left|Z(\sigma)\right|\right|_{L^4}^3\right)d\sigma\right)^{\frac{2}{3}} \times \|w\|_{L^4} 
\end{align*}
and use $Z^{n(k)} \to Z$ in $L^3(0,T,L^4)$. 
Therefore, we conclude, as $k\to \infty,$
\begin{align}
\left(X(t)-u(0), w\right)=\int_0^t\left(X(\sigma),\mathcal{L}_{1}w\right)d\sigma-\int_0^t\left((X(\sigma)+Z(\sigma))^3,w\right)d\sigma.\,\,\,\mathbb{P}'-a.s.
\end{align}

We have proved the existence of the solution on $(\Omega',\mathcal{F}',\mathbb{P}').$
To obtain the solution on the original space $(\Omega,\mathcal{F},\mathbb{P}),$ we use the Gy\"ongy-Krylov criterion. 
Regarding $X_{n}$ of Proposition \ref{prop:gk} as $(v_{n(k)},\Pi_{n(k)}Z),$ we know that arbitrary subsequence of  $(v_{n(k)},\Pi_{n(k)}Z)$ converges to $\mu$ in law. 
It is thus sufficient to prove that for any $\varepsilon>0$ (see the details in Section 4.4.2 of \cite{FFL}),
\begin{align*}
\lim_{h\rightarrow\infty}\mathbb{P}(\|(v_{n_1(k(h))},\Pi_{n_1(k(h))})-(v_{n_2(k(h))},\Pi_{n_2(k(h))}))\|_{O}>\varepsilon)=0.
\end{align*}
with 
$$O:= L^3(0,T;L^4)\cap C([0,T]; \mathcal{W}^{-(1+\alpha),2}) \cap L^2(0,T;\mathcal{W}^{1-\alpha,2}) \times C([0,T];L^4).$$
This follows from the equivalence of probability laws between $(X^{n(k)},Z^{n(k)})$ and $(v_{n(k)},\Pi_{n(k)}Z),$ and convergence 
\begin{align*}
&\lim_{k\rightarrow\infty}(X^{n(k)},Z^{n(k)})=(X,Z), \quad {\rm in}\\ 
&L^3(0,T;L^4)\cap C([0,T]; \mathcal{W}^{-(1+\alpha),2}) \cap L^2(0,T;\mathcal{W}^{1-\alpha,2} ) \times C(0,T;L^4),\,\,\mathbb{P}'-a.s.
\end{align*} 
Therefore there exists a random variable $(V,Z)$ such that $(v_{n(k)},\Pi_{n(k)}Z)\rightarrow (V,Z)$ in $L^3(0,T;L^4) \cap C([0,T]; \mathcal{W}^{-(1+\alpha),2}) \cap L^2(0,T;\mathcal{W}^{1-\alpha,2}) \times C(0,T;L^4)$ in probability. Note that by taking a subsequence $\{v_{n(k(l))}\}$, the convergence in probability becomes $\mathbb{P}$-a.s. convergence. 
Similarly to the above discussion, we get

\begin{align}
\left(V(t)-u(0), w\right)=\int_0^t\left(V(u),\mathcal{L}_{1}w\right)du-\int_0^t\left((V(u)+Z(u))^3,w\right)du.\,\,\,\mathbb{P}-a.s.
\end{align}
Thus we have proved the existence of the solution on $(\Omega, \mathcal{F}, \prob)$.

Recall that we proved $v_{n(k)}\rightarrow V$ in $L^3(0,T;L^4)\cap C([0,T]; \mathcal{W}^{-(1+\alpha),2}) \cap L^2(0,T,\mathcal{W}^{1-\alpha,2})$, $\mathbb{P}$ almost surely. Also, 
\begin{align}\label{BDD}
\sup_{n\in\mathbb{N}}\mathbb{E}\left(\sup_{t\in[0,T]}||v_n(t)||_{L^2}^2\right)\le C(T,||u_0||_{L^2}).
\end{align} 
This inequality comes from $(\ref{EE2n})$.The inequality 
$(\ref{BDD})$ implies that $\{v_{n(k)}\}$ is weak star compact in $L^2(\Omega, L^{\infty}(0,T;L^2))$. 
\textcolor{black}{Thus there exist a subsequence (denoted by the same letter) and a limit $v^*$ such that 
$v_{n(k)} \to v^*$ weak star in $L^2(\Omega; L^{\infty}(0,T;L^2))$.} 
On the other hand, $v_{n(k(l))}\rightarrow V$ in $L^3(0,T;L^4)\cap C([0,T]; \mathcal{W}^{-(1+\alpha),2}) \cap
L^2(0,T, \mathcal{W}^{1-\alpha,2})$ $\mathbb{P}$-almost surely. 
In particular, $v_{n(k(l))}\rightarrow V$ weak star in $L^3(0,T;L^4) 
\cap C([0,T]; \mathcal{W}^{-(1+\alpha),2}) \cap L^2(0,T;\mathcal{W}^{1-\alpha,2})$ almost surely.
We know 
\begin{align*}
v^*\in L^{\infty}(0,T;L^2)\subset L^3(0,T;L^4)\cap C([0,T]; \mathcal{W}^{-(1+\alpha),2}) 
\cap L^2 (0,T;\mathcal{W}^{1-\alpha,2}),
\end{align*}
almost surely.
The uniquness of weak star limit implies
\begin{align*}
v^*=V.
\end{align*}
Therefore,
\begin{align*}
V\in L^3(0,T;L^4)\cap C([0,T]; \mathcal{W}^{-(1+\alpha),2})\cap L^2(0,T, \mathcal{W}^{1-\alpha,2}) 
\cap L^{\infty}(0,T:L^2),
\end{align*}
$\prob$- almost surely. \textcolor{black}{Finally we show the pathwise uniqueness of solutions following the idea in \cite{BBS}. Consider two solution $v_1, v_2 \in L^{\infty}(0,T;L^2)$, and set $d=v_1-v_2$. Then $d$ satisfies
$$ \partial_t d = {\mathcal{L}}_1 d -\{(v_1+Z)^3 -(v_2+Z)^3\}, $$
and
\begin{equation} \label{d}
\frac12 \partial_t \|d\|_{L^2}^2 = (\partial_t d,d)=(\mathcal{L}_1 d,d)-\int_{\mathbb{T}^2} d \{d^3 +3d^2(v_2+Z) +3d(v_2+Z)^2\}dx.
\end{equation}
Note that 
$$-3 \int_{\mathbb{T}^2} d^3 (v_2+Z) dx \le \int_{\mathbb{T}^2} \left(d^4 +\frac94 d^2 (v_2+Z)^2 \right)dx, $$
and 
\begin{eqnarray*}
 (\mathcal{L}_1 d,d) &=& 2 \|d\|_{L^2}^2-\|(1-\mathcal{L}_0)^{\frac12} d\|_{L^2}^2 \\
 &=& 2 \|d\|_{L^2}^2 -\|d\|_{\mathcal{W}^{1,2}}^2.
\end{eqnarray*}
Thus, 
$$ (\ref{d}) \le \{2+\frac94 \sup_{t\in [0,T]} (\|v_2\|_{L^2}^2 +\|Z\|_{L^2}^2 )\} \|d\|_{L^2}^2, $$
which implies $d=0$ in $L^2$ after an application of the Gronwall inequality. 
}
This completes the proof of Theorem $\ref{V}$.$\qed$
%
%

\section{Numerical simulation}
In this section, we present some simulations in space dimension 2. The idea is to first perform simulations for the equation of $A$ and convert them to $u$ by the Ansatz, and at the same time, to perform direct simulations for the Swift-Hohenberg equation on $u$. We expect that with the proposed definitions of the noise term and the scaling, the patterns obtained by the two methods will be similar one to the other. And as informal observations, we perform simulations to compare the results in both the deterministic and stochastic case.

\subsection*{Equation of $A$}
We recall the equation of $A(X,Y,T)$ in the deterministic case
\begin{equation}\label{equation_A}
	\partial_T A =4\partial^2_X A +\partial_Y^2 A +A - 3\lvert A\rvert^2A,
\end{equation}	
in the space domain $[-L,L]\times[-L,L]$ and time interval $[0,T_A]$. As $A$ is complex-valued, we suppose $A ={A^R}+\mathbf{i}{A^I} $ , where $R$ stands for real and $I$ stands for imaginary. And we separate the real part and the imaginary part of the equation into the following system
\begin{equation}\label{equation_A_real_imaginary}
	\begin{dcases}
		\partial_T {A^R}=4\partial^2_X {A^R}+\partial_Y^2 {A^R}+{A^R}- 3\lvert ({A^R})^2 + ({A^I})^2\rvert {A^R},\\
		\partial_T {A^I} =4\partial^2_X {A^I} +\partial_Y^2 {A^I} +{A^I} - 3\lvert ({A^R})^2 + ({A^I})^2\rvert {A^I},
	\end{dcases}
\end{equation}	
with initial conditions and periodic boundary conditions for both ${A^R}$ and ${A^I}$.
\subsubsection*{Convert $A$ to $u$ by the Ansatz}
Once we obtain the numerical solution of $A$, we convert it to $u$ by applying the following Ansatz
\begin{equation}\label{Ansatz_A_to_U}
	u(x,y,t) = \varepsilon A(\varepsilon x,\varepsilon y, \varepsilon^2 t)e^{\mathbf{i}x}+\varepsilon \overline{A}(\varepsilon x,\varepsilon y, \varepsilon^2 t) e^{-\mathbf{i}x}
\end{equation}
where $x=X/\varepsilon $,  $y=Y/\varepsilon $, $t = T/\varepsilon^2$ and $ \overline{A}={A^R}-\mathbf{i}{A^I} $.
A direct computation yields
\begin{equation}\label{Ansatz_A_to_U_detail}
	u(x,y,t) = 2\varepsilon \left(  {A^R}(\varepsilon x,\varepsilon y, \varepsilon^2 t)\cdot \cos(x) - {A^I}(\varepsilon x,\varepsilon y, \varepsilon^2 t)\cdot \sin(x)\right).
\end{equation}
\subsection*{Anisotropic Swift-Hohenberg equation}
The anisotropic Swift-Hohenberg equation that we consider is as follows
\begin{equation}\label{2D_SH_anisotropic}
	\partial_t u = -(1+\partial_x^2)^2 u +	\partial_y^2 u + \varepsilon^2 u -u^3.
\end{equation}
In order to perform simulations, we decompose the equation into the system
\begin{equation}\label{2D_SH_decomposite}
	\begin{dcases}
		\partial_t u&= -u+ \partial_x^2 \mu  +\partial^2_{y}u+\varepsilon^2u  -  u^3\\
		\mu &= -\partial_x^2 u - 2 u
	\end{dcases}
\end{equation}
with initial condition for $u$ and periodic boundary condition for both $u$ and $\mu$.
\subsection*{Form of the stochastic term}
In order to consider the corresponding stochastic equations of \eqref{equation_A} and \eqref{2D_SH_anisotropic}, we present the stochastic term. We first define the stochastic term for $A$, which is
\begin{equation}\label{stochastic_term_for_A}
	\xi(X, Y, T) =C_L \sum_{\mathbf{k}\in \mathbb{Z}^2} \beta_{\mathbf{k}}(T)e^{\frac{\mathbf{-i \pi \mathbf{k}}}{L} \cdot \mathbf{X}}.
\end{equation}	
More precisely, we have the space domain $[-L,L]\times [-L,L]$ and $C_L=1/(2L)$, $\mathbf{k}={k^R}+\mathbf{i} {k^I}$, $\mathbf{X}=(X, Y)$ and $\beta_{\mathbf{k}}= \beta^R_{{({k^R},{k^I})}}+ \mathbf{i}\beta^I_{{({k^R},{k^I})}}$, where  $\beta^R_{{({k^R},{k^I})}}$ and $ \beta^I_{{({k^R},{k^I})}}$ are independent real-valued Brownian motions. 
The corresponding stochastic equation of \eqref{equation_A} is given by
\begin{equation}\label{equation_A_sto}
	\partial_T A =4\partial^2_X A +\partial_Y^2 A +A - 3\lvert A\rvert^2A+\dot{\xi},
\end{equation}
with 
\begin{equation}\label{noise_term_for_A}
	\dot{\xi}(X, Y, T) = C_L \sum_{\mathbf{k}\in \mathbb{Z}^2}\dot{\beta_{\mathbf{k}}}(T)e^{\frac{\mathbf{-i \pi \mathbf{k}}}{L} \cdot \mathbf{X}}.
\end{equation}	

Since $\xi(X, Y, T)$ is complex-valued, we suppose $\xi(X, Y, T)=\xi^{R}(X, Y, T)+\mathbf{i}\xi^{I}(X, Y, T)$. A computation yields
\begin{equation}
	\begin{dcases}
		\xi^{R} = \sum_{\substack{{{k^R}}\in \mathbb{Z}\\{{k^I}}\in \mathbb{Z}}} \left[\beta^R_{{({k^R},{k^I})}}(T)
		\cos\Big(\frac{{ \pi ({k^R}X+{k^I}Y )}}{L} \Big) -
		\beta^I_{{({k^R},{k^I})}}(T)\sin\Big(\frac{{\pi  ({k^R}X+{k^I}Y)}}{L}  \Big)  \right]\\
		\xi^{I} = \sum_{\substack{{{k^R}}\in \mathbb{Z}\\{{k^I}}\in \mathbb{Z}}}\left[\beta^R_{{({k^R},{k^I})}}(T) \sin\Big(\frac{  \pi ({k^R}X+{k^I}Y)}{L} \Big) +
		\beta^I_{{({k^R},{k^I})}}(T)\cos\Big(\frac{ \pi ({k^R}X+{k^I}Y}{L}\Big)  \right].
	\end{dcases}
\end{equation}
And the corresponding stochastic system of \eqref{equation_A_real_imaginary} is given by
\begin{equation}\label{equation_A_real_imaginary_sto}
	{\displaystyle \begin{cases}
			\partial_T {A^R}=4\partial^2_X {A^R}+\partial_Y^2 {A^R}+{A^R}- 3\lvert ({A^R})^2 + ({A^I})^2\rvert {A^R} + \dot{\xi}^{R}\\
			\partial_T {A^I} =4\partial^2_X {A^I} +\partial_Y^2 {A^I} +{A^I} - 3\lvert ({A^R})^2 + ({A^I})^2\rvert {A^I} +\dot{\xi}^{I}.
	\end{cases}}
\end{equation}	

Next, we present the stochastic term for the equation of $u$,  after some computation, we define
\begin{equation}\label{noise_form_for_u}
	\xi_{\varepsilon}(x,y,t) =  C_L \sum_{\mathbf{k}\in \mathbb{Z}^2} \beta_{\mathbf{k}}(\varepsilon^2 t)e^{\frac{\mathbf{-i \pi \mathbf{k}}}{L/\varepsilon} \cdot \mathbf{x}},
\end{equation}	
where $\mathbf{x}=(x,y) \in [-L/\varepsilon,L/\varepsilon]\times[-L/\varepsilon,L/\varepsilon]$ with $x=X/\varepsilon $,  $y=Y/\varepsilon $ and $t = T/\varepsilon^2$.
We suppose that $ {\beta}_\mathbf{-k}= \overline{\beta}_\mathbf{k}$ and as a result, 
\begin{equation}\label{noise_form_for_u_detailed}
	\begin{aligned}
		&	 \xi_{\varepsilon}(x,y,t) \\= & 2  C_L \sum_{\substack{{k^R} \in \mathbb{Z}\\{k^I} \in \mathbb{Z}^+}}
		\left(
		\beta^R_{{({k^R},{k^I})}}(\varepsilon^2 t)\cos\Big( \frac{\pi ({k^R} x+{k^I} y)}{L/\varepsilon}  \Big) +\beta^I_{{({k^R},{k^I})}}(\varepsilon^2 t)\sin\Big(\frac{\pi ({k^R} x+{k^I} y)}{L/\varepsilon}\Big)  
		\right)
	\end{aligned}
\end{equation}
and we remark that $\varepsilon\xi_{\varepsilon}$ is real-valued. The stochastic system corresponds to \eqref{2D_SH_decomposite} is given by
\begin{equation}\label{2D_SH_decomposite_sto}
	\begin{dcases}
		\partial_t u&= -u+ \partial_x^2 \mu  +\partial^2_{y}u+\varepsilon^2u  -  u^3 +\varepsilon \cdot \dot{\xi}_{\varepsilon}\\
		\mu &= -\partial_x^2 u - 2 u.
	\end{dcases}
\end{equation}


\subsection{Space and time discretizations}
We mainly present the numerical settings of $A$ and the settings for $u$ are defined correspondingly. We discretize the space domain $[-L,L]\times[-L,L]$ into a $N_X \times N_Y$ uniform mesh, so that $\Delta X= 2L/N_X$ and $\Delta Y= 2L/N_Y$. 
We define $p\in \{0,1,2,...,N_X-1\}$ and $q \in \{0,1,2,...,N_Y-1\}$, two indices in direction $X$ and $Y$ respectively. The control volume $(p,q)$ is the volume whose barycenter satisfies $$\mathbf{X}_{p,q}=\left((p+0.5)\cdot \Delta X, \linebreak (q+0.5)\cdot \Delta Y\right). $$

And we apply uniform time discretization, that is we fix the time step $\Delta T$ and define $T_n =n\Delta T$ for all $n=0,1,2,...$. If we consider $N_T$ time steps, the total time interval is $ \cup_{n=0}^{N_T-1}[n\Delta T,(n+1)\Delta T)$.

The discrete solutions of ${A^R}$ and of ${A^I}$ are denoted by $\{{A^R}_{p,q}^n\}$ and $\{{A^I}_{p,q}^n\}$ over control volume ${(p,q)}$ in time interval $[T_n, T_{n+1})$ respectively.

In the discretization of the space derivative, we will need the values of ${A^R}_{-1,q}^n$, ${A^R}_{N_X,q}^n$, ${A^R}_{p,-1}^n$ and ${A^R}_{p,N_Y}^n$, because of the periodic boundary condition, we set
\begin{equation}\label{discrete_periodic_boundary_x}
	\begin{aligned}
		{A^R}_{-1,q}^n:={A^R}_{N_X-1,q}^n \quad\mbox{ for all }\quad q \in \{0,1,2,...,N_Y-1\}\\
		{A^R}_{N_X,q}^n:=	{A^R}_{0,q}^n \quad\mbox{ for all }\quad q \in \{0,1,2,...,N_Y-1\}\\
	\end{aligned}
\end{equation}
and
\begin{equation}\label{discrete_periodic_boundary_y}
	\begin{aligned}
		{A^R}_{p,-1}^n:={A^R}_{p,N_Y-1}^n \quad\mbox{ for all }\quad p\in \{0,1,2,...,N_X-1\},\\
		{A^R}_{p,N_Y}^n:={A^R}_{p,0}^n \quad\mbox{ for all }\quad p\in \{0,1,2,...,N_X-1\},
	\end{aligned}
\end{equation}
and we have the same conditions for $\{ {A^I}_{p,q}^n\}$.
For the approximation of $u$, we apply corresponding settings with notations $\Delta x$, $\Delta y$ and $\Delta t$, $N_t$ and the discrete solution is denoted by $\{u_{p,q}^n\}$.
\subsection{Discretization of the noise term}
Suppose $\beta$ is a Brownian motion,
for the numerical simulations, we approximate $\dot{\beta}$ by
\begin{equation*}
	\dot{\beta}(t) \approx \frac{\beta(t+\Delta t)-\beta(t)}{\Delta t},
\end{equation*}
where $\beta(t+\Delta t)-\beta(t) \sim \mathcal{N}(0,\Delta t)$ is a Gaussian random variable with mean value $0$ and variance $\Delta t$.

We discretize the noise term \eqref{noise_term_for_A} as follows
\begin{equation}
	\dot{\xi}(X,Y,T)\\ \approx  C_L\sum_ { \substack{{k^R}\in \{-m_{\mathrm{R}},  ...,0,...,m_{\mathrm{R}}\} \\ {{k^I}\in \{-m_{\mathrm{I}},  ...,0,...,m_{\mathrm{I}}\} }}}
	\frac{\beta_{\mathbf{k}}(T+\Delta T) -  \beta_{\mathbf{k}}(T) }{\Delta T}e^{\frac{\mathbf{-i \pi \mathbf{k}}}{L} \cdot \mathbf{X}}
\end{equation}
such that
\begin{equation}
	\begin{aligned}
		&\dot{\xi^{R}}(X,Y,T)\\ \approx & C_L\sum_ { \substack{{k^R}\in \{-m_{\mathrm{R}},  ...,0,...,m_{\mathrm{R}}\} \\ {{k^I}\in \{-m_{\mathrm{I}},  ...,0,...,m_{\mathrm{I}}\} }}}
		\Bigg[ \frac{\beta^R_{({k^R},{k^I})}(T+\Delta T) -  \beta^R_{({k^R},{k^I})}(T) }{\Delta T}
		\cos\Big(\frac{{ \pi ({k^R}X+{k^I}Y )}}{L} \Big) \\&\quad\quad \quad\quad\quad\quad\quad\quad\quad\quad -
		\frac{\beta^I_{({k^R},{k^I})}(T+\Delta T) -  \beta^I_{({k^R},{k^I})}(T) }{\Delta T}\sin\Big(\frac{{\pi  ({k^R}X+{k^I}Y)}}{L}  \Big)  \Bigg]\\
	\end{aligned}
\end{equation}
and
\begin{equation}
	\begin{aligned}
		&\dot{\xi^{I}}(X,Y,T)\\ \approx & C_L\sum_ { \substack{{k^R}\in \{-m_{\mathrm{R}},  ...,0,...,m_{\mathrm{R}}\} \\ {{k^I}\in \{-m_{\mathrm{I}},  ...,0,...,m_{\mathrm{I}}\} }}}
		\Bigg[\frac{\beta^R_{({k^R},{k^I})}(T+\Delta T) -  \beta^R_{({k^R},{k^I})}(T) }{\Delta T}
		\sin\Big(\frac{{ \pi ({k^R}X+{k^I}Y )}}{L}\Big) \\&\quad\quad\quad\quad\quad\quad\quad\quad \quad\quad +
		\frac{\beta^I_{({k^R},{k^I})}(T+\Delta T) -  \beta^I_{({k^R},{k^I})}(T) }{\Delta T}\cos\Big(\frac{{\pi  ({k^R}X+{k^I}Y)}}{L}  \Big)  \Bigg]\\
	\end{aligned}
\end{equation}
where $\beta^R_{({k^R},{k^I})}(T+\Delta T) -  \beta^R_{({k^R},{k^I})}(T)$ and $\beta^I_{({k^R},{k^I})}(T+\Delta T) -  \beta^I_{({k^R},{k^I})} \sim \mathcal{N}(0, \Delta T)$. The $m_{\mathrm{R}}$ and $m_{\mathrm{I}}$ are the truncation numbers. We denote this approximation of $\dot{\xi}(X,Y,T)$ by $\Xi$ such that $\dot{\xi^{R}}(X,Y,T)$ is approximated by $\Xi^{R}$ and $\dot{\xi^{I}}(X,Y,T)$ by $\Xi^{I}$.
%
In view of \eqref{noise_form_for_u_detailed}, the discretization of the noise term \eqref{noise_form_for_u} is as follows
\begin{equation}
	\begin{aligned}
		&\dot{\xi}_{\varepsilon}(x,y,t)\\ \approx & 2C_L\sum_ { \substack{{k^R}\in \{-m_{\mathrm{R}},  ...,0,...,m_{\mathrm{R}}\} \\ {{k^I}\in \{0,...,m_{\mathrm{I}}\} }}}
		\Bigg[\frac{\beta^R_{({k^R},{k^I})} (\varepsilon^2(t+\Delta t) )-  \beta^R_{({k^R},{k^I})}(\varepsilon^2 t) }{\Delta t}
		\cos\Big(\frac{{ \pi ({k^R} x+{k^I} y )}}{L/\varepsilon} \Big) \\&\quad\quad\quad\quad\quad\quad\quad\quad\quad\quad +
		\frac{\beta^I_{({k^R},{k^I})} (\varepsilon^2(t+\Delta t) )-  \beta^I_{({k^R},{k^I})}(\varepsilon^2 t) }{\Delta t}\sin\Big(\frac{{\pi  ({k^R}x+{k^I}y)}}{L/\varepsilon}  \Big)  \Bigg]\\
	\end{aligned}
\end{equation}
where $\beta_{{k^R}}(\varepsilon^2 (t+\Delta t)) -  \beta_{{k^R} }(\varepsilon^2 t)$ and $\beta_{{k^I}}(\varepsilon^2 (t+\Delta t)) -  \beta_{{k^I}}(\varepsilon^2 t)  \sim \mathcal{N}(0, \varepsilon^2\Delta t)$. We denote this approximation term by $\Xi_{\varepsilon}$.

We present in the following directly the numerical scheme for the stochastic case. 
For the deterministic case, we perform the simulation by omitting the stochastic term.
\subsection{Numerical schemes}
We apply the finite difference scheme and first present the scheme for $A$ of the system \eqref{equation_A_real_imaginary_sto}
\begin{equation}\label{scheme_A_real_imaginary_sto}
	\begin{dcases}
		\frac{ {A^R}^{n+1}_{p,q}-{A^R}^{n}_{p,q}}{\Delta T} =4d^2_X {A^R}^{n+1}_{p,q}
		+d^2_Y {A^R}^{n}_{p,q}+{A^R}^{n+1}_{p,q}- 3\lvert ({A^R}^{n}_{p,q})^2 + ({A^I}^{n}_{p,q})^2\rvert {A^R}^{n}_{p,q}  + \Xi^{R}\\
		\frac{ {A^I}^{n+1}_{p,q}-{A^I}^{n}_{p,q}}{\Delta T} =4d^2_X {A^I}^{n+1}_{p,q}
		+d^2_Y {A^I}^{n}_{p,q}+{A^I}^{n+1}_{p,q}- 3\lvert ({A^R}^{n}_{p,q})^2 + ({A^I}^{n}_{p,q})^2\rvert {A^I}^{n}_{p,q}  + \Xi^{I},
	\end{dcases}
\end{equation}	
where $d^2_X$ and $d^2_Y$ are discrete operators such that 
$$d^2_X {A^R}^{n+1}_{p,q}=\frac{{A^R}^{n+1}_{p-1,q}+{A^R}^{n+1}_{p+1,q}- 2{A^R}^{n+1}_{p,q}   }{(\Delta X)^2} \quad \mbox{and} \quad
d^2_Y {A^R}^{n}_{p,q} =\frac{{A^R}^{n}_{p,q-1}+{A^R}^{n}_{p,q+1}- 2{A^R}^{n}_{p,q}  }{(\Delta Y)^2},
$$
for all $p\in \{0,1,2,...,N_X-1\}$ and $q \in \{0,1,2,...,N_Y-1\}$. We refer to \eqref{discrete_periodic_boundary_x} and \eqref{discrete_periodic_boundary_y} for the periodic boundary condition. For $n=0,1,2,...,N_T-1$, knowing the values of $({A^R}^{n}_{p,q}, {A^I}^{n}_{p,q})$, we compute the values of $({A^R}^{n+1}_{p,q},{A^I}^{n+1}_{p,q})$, for all $p\in \{0,1,2,...,N_X-1\}$ and $q \in \{0,1,2,...,N_Y-1\}$.

And we implement the following numerical scheme for the Swift-Hohenberg equation in the form of system \eqref{2D_SH_decomposite_sto}.
\begin{equation}\label{scheme_u}
	\begin{dcases}
		\frac{u_{p,q}^{n+1}-u_{p,q}^{n}}{\Delta t} &=  -u_{p,q}^{n+1} +d^2_x \mu^{n+1}_{p,q} + d^2_y u^{n}_{p,q} +\varepsilon^2 u^n_{p,q}- (u^n_{p,q})^3 +\varepsilon\Xi_{\varepsilon}\\
		\mu_{p,q}^{n+1} & =-d^2_x u^{n+1}_{p,q} -  u^{n+1}_{p,q}- u^n_{p,q},
	\end{dcases}
\end{equation}
with corresponding definitions of $d^2_x$ and $d^2_y$ and the periodic boundary conditions.
For $n=0,1,2,...,N_t-1$, knowing the values of $(u^{n}_{p,q}, \mu^{n}_{p,q})$, we compute the values of $(u^{n+1}_{p,q},\mu^{n+1}_{p,q})$, for all $p\in \{0,1,2,...,N_x-1\}$ and $q \in \{0,1,2,...,N_y-1\}$.


\subsection{Numerical settings}
In the simulations, $\varepsilon$ is the parameter to connect $A$ and $u$, so we first fix the value of $\varepsilon$. 
And we perform simulations for $A$ with the following settings.\\
\textit{Numerical settings for $A$}
\begin{itemize}
	\item The space domain to be $ [-L, L]\times[-L, L]$ with $L=\pi/2$;
	\item we discretize the space into $100\times 100$ uniform square;
	\item we fix the time step $\Delta t = 0.0001$;
	\item we set the initial condition ${A^R}(X,Y,0)=A^R_{0}(X,Y)$ and ${A^I}=A^I_{0}(X,Y)$;
	\item we perform simulations of $A$ and convert the numerical results to $u$ by the Ansatz.
\end{itemize}
We perform simulations of $u$ by \eqref{scheme_u} with the following settings.\\
\textit{Numerical settings for $u$}
\begin{itemize}
	\item The space domain to be $ [-L/\varepsilon, L/\varepsilon]\times[-L/\varepsilon, L/\varepsilon]$ with $L=\pi/2$;
	\item we discretize the space into $100\times 100$ uniform square;
	\item we choose time step $\Delta t = 0.001$;
	\item we compute the initial condition for $u$ based on the initial condition of ${A^R}$ and ${A^I}$ by the Ansatz \eqref{Ansatz_A_to_U}, which yields
	\begin{equation}\label{Ansatz_A_to_U_initial_condition}
		u_0(x,y) = 2\epsilon \left(  {A^R}(\varepsilon x,\varepsilon y, 0)\cdot \cos(x) - {A^I}(\varepsilon x,\varepsilon y, 0)\cdot \sin(x)\right);
	\end{equation}
	\item we perform simulations for $u$.
\end{itemize}	
\subsection{Results and observations}

We set  $\varepsilon=0.25$ and perform numerical simulations for ${A^I}$ and ${A^R}$ with the initial conditions
\begin{equation}\label{initial_condition_{A^R}}
	{A^R_0}(X,Y) =
	\begin{dcases}
		1 \quad& \mbox{if}\quad Y\in (-L,0)\\
		0 \quad &\mbox{otherwise}\quad
	\end{dcases}
\end{equation}
and
\begin{equation}\label{initial_condition_{A^I}}
	{A^I_0}(X,Y) =
	\begin{dcases}
		0 \quad& \mbox{if}\quad Y\in (-L,0)\\
		1 \quad &\mbox{otherwise.}\quad
	\end{dcases}
\end{equation}
The result is as follows.
\subsubsection*{Deterministic case}		
\begin{figure}[H]
	\centering
	\includegraphics[width=12cm]{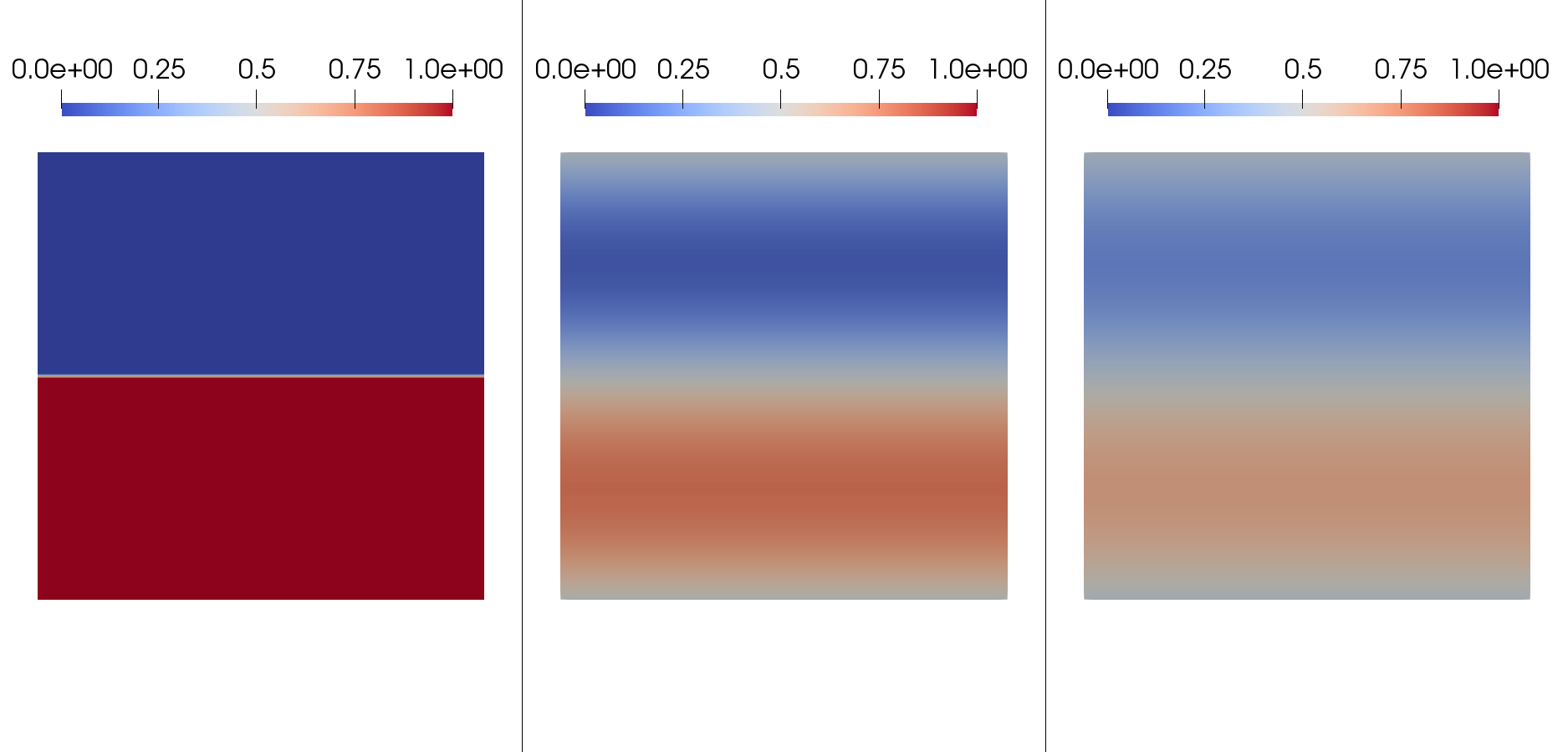}
	\caption{Solution for ${A^R}$ at time $T=0$, $0.1$ and $0.2$}
\end{figure}
\begin{figure}[H]
	\centering
	\includegraphics[width=12cm]{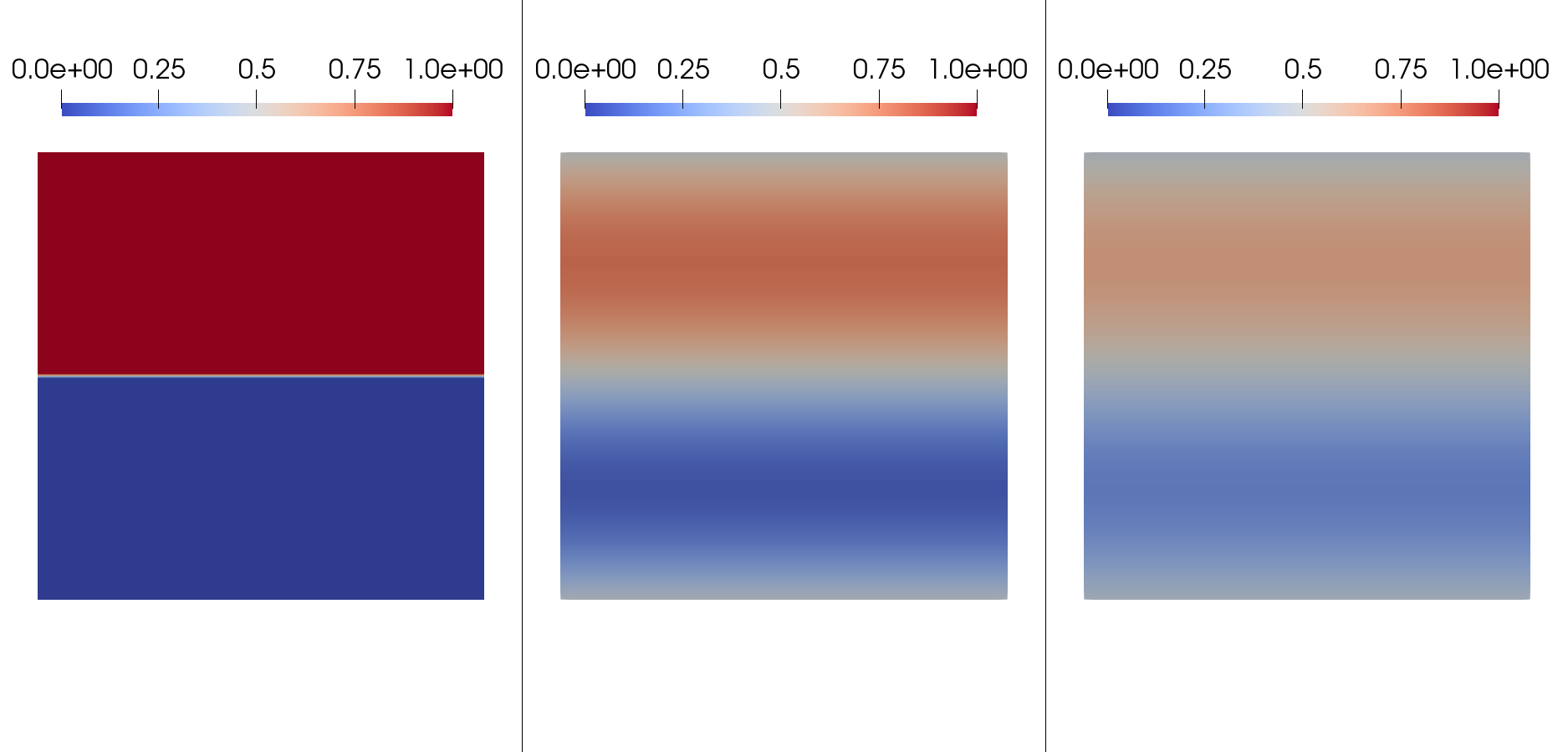}
	\caption{Solution for ${A^I}$ at time $T=0$, $0.1$ and $0.2$}
\end{figure}
We convert the numerical results from $A$ to $u$.
\begin{figure}[H]
	\centering
	\includegraphics[width=12cm]{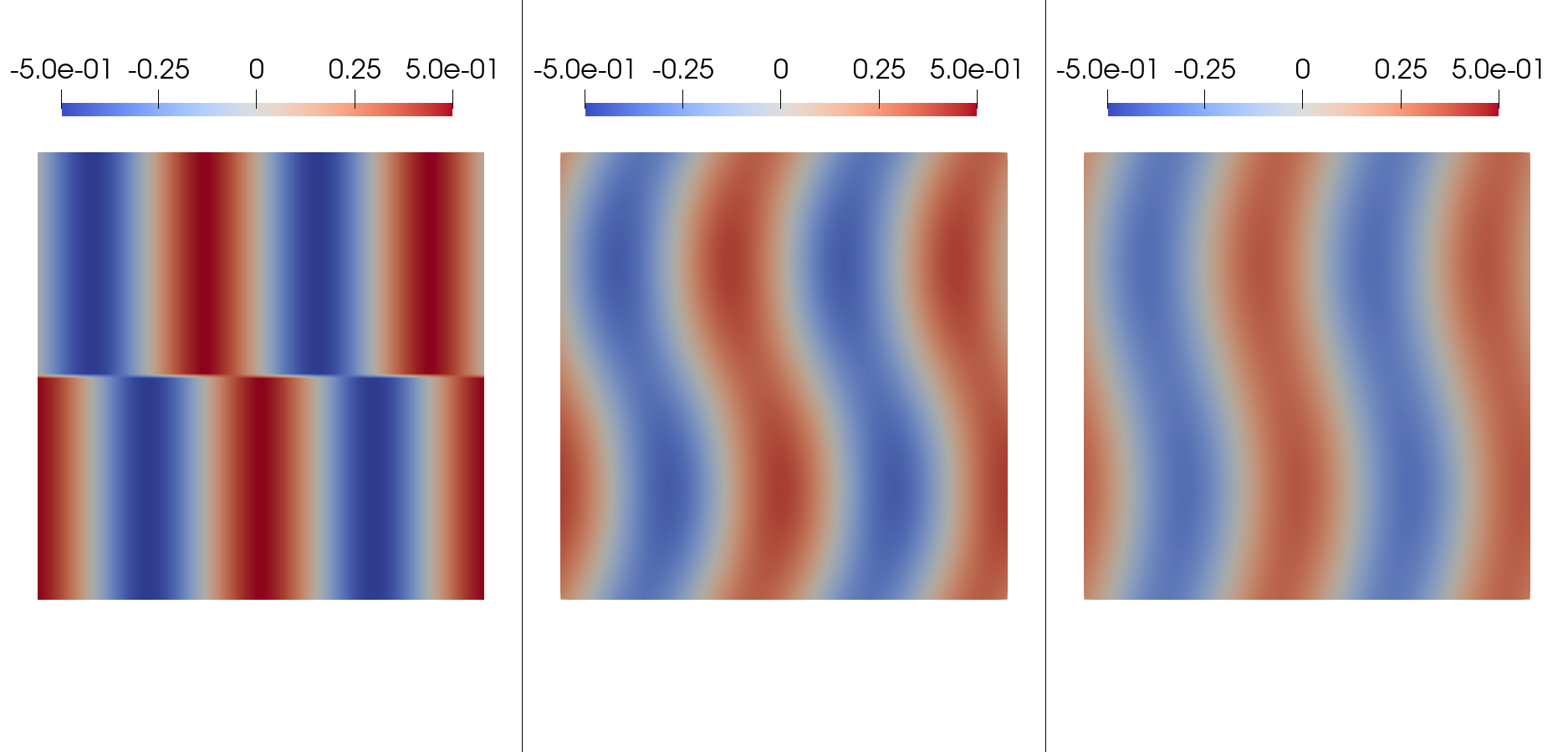}
	\caption{$u$ at time $t=0,0.1/\varepsilon^2$ and $0.2/\varepsilon^2$ converted by the Ansatz from $A$ of time $T=0$, $0.1$ and $0.2$}
	\label{u_byAnsatz_det_eps005}
\end{figure}
Then we perform directly the simulation for $u$ and compare it to the simulation results of $u$ convert by $A$ which are presented in Figure \ref{u_byAnsatz_det_eps005}.
\begin{figure}[H]
	\centering
	\includegraphics[width=12cm]{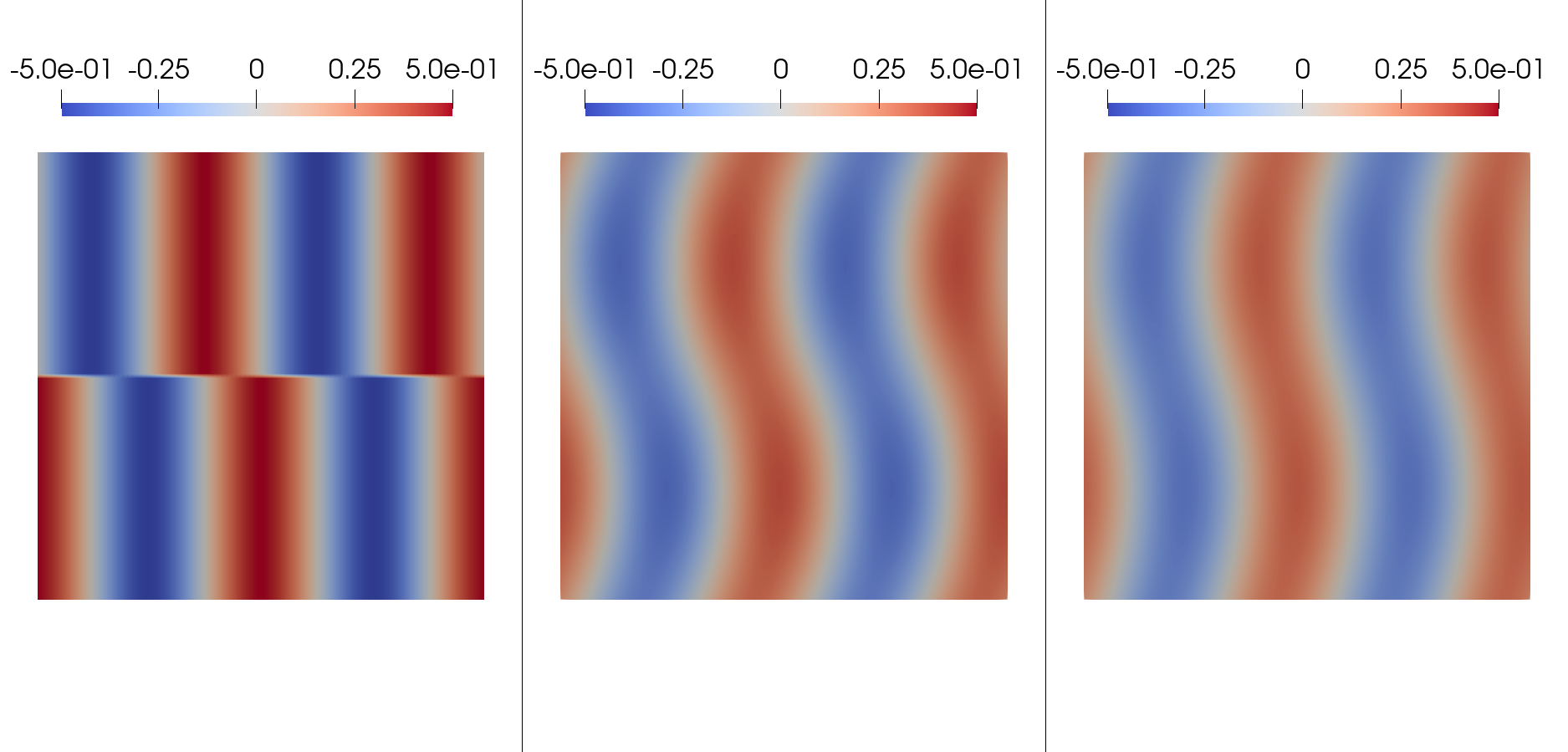}
	\caption{$u$ by direct simulation at time $t=0, 0.1/\varepsilon^2$ and $0.2/\varepsilon^2$}
\end{figure}

\subsubsection*{Stochastic case}
\textcolor{black}{We set the truncation numbers $m_{\mathrm{R}}=10$ and $m_{\mathrm{I}}=10$.}
We first present the simulations results of ${A^R}$ and 
${A^I}$
\begin{figure}[H]
	\centering
	\includegraphics[width=12cm]{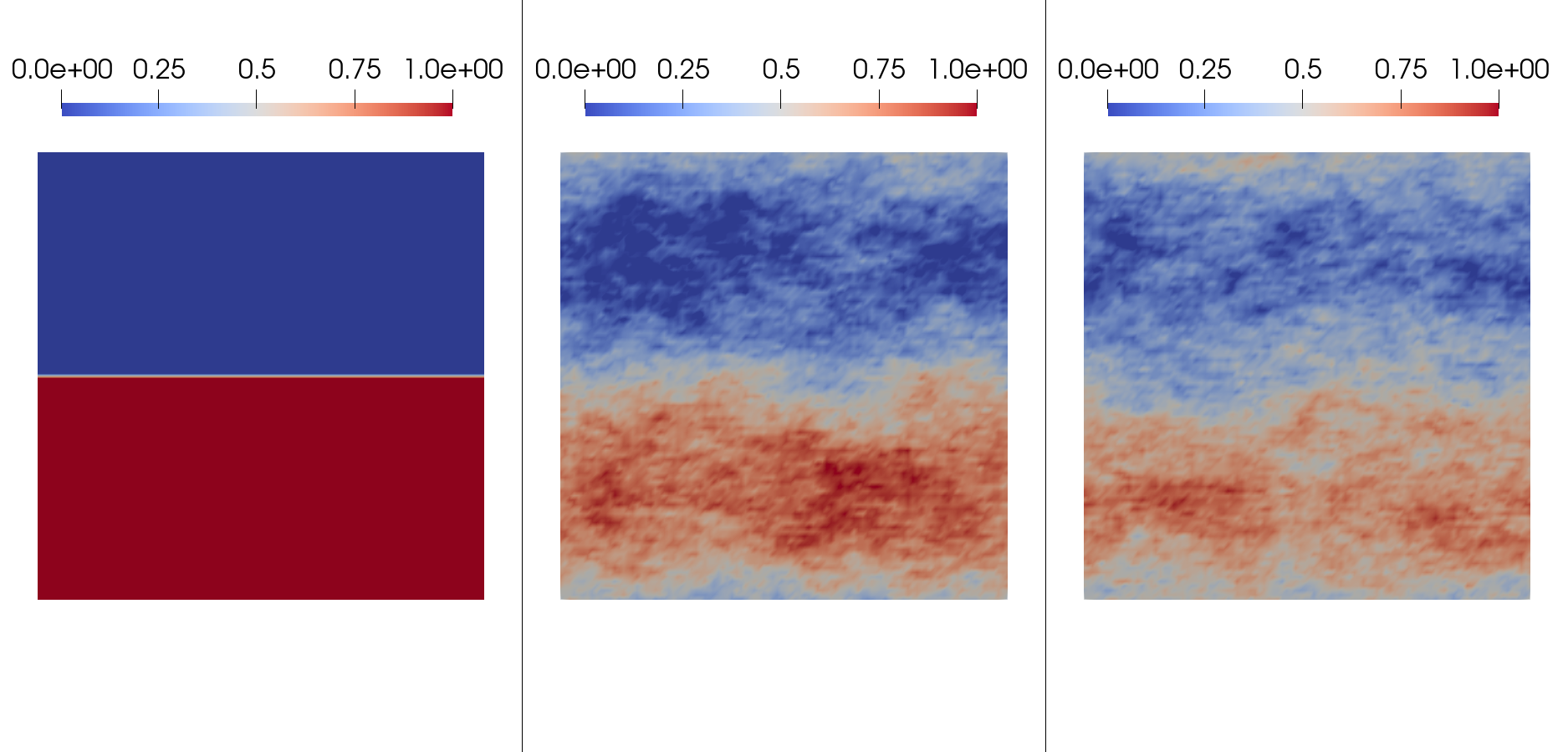}
	\caption{Solution for ${A^R}$ at time $T=0$, $0.1$ and $0.2$}
\end{figure}
\begin{figure}[H]
	\centering
	\includegraphics[width=12cm]{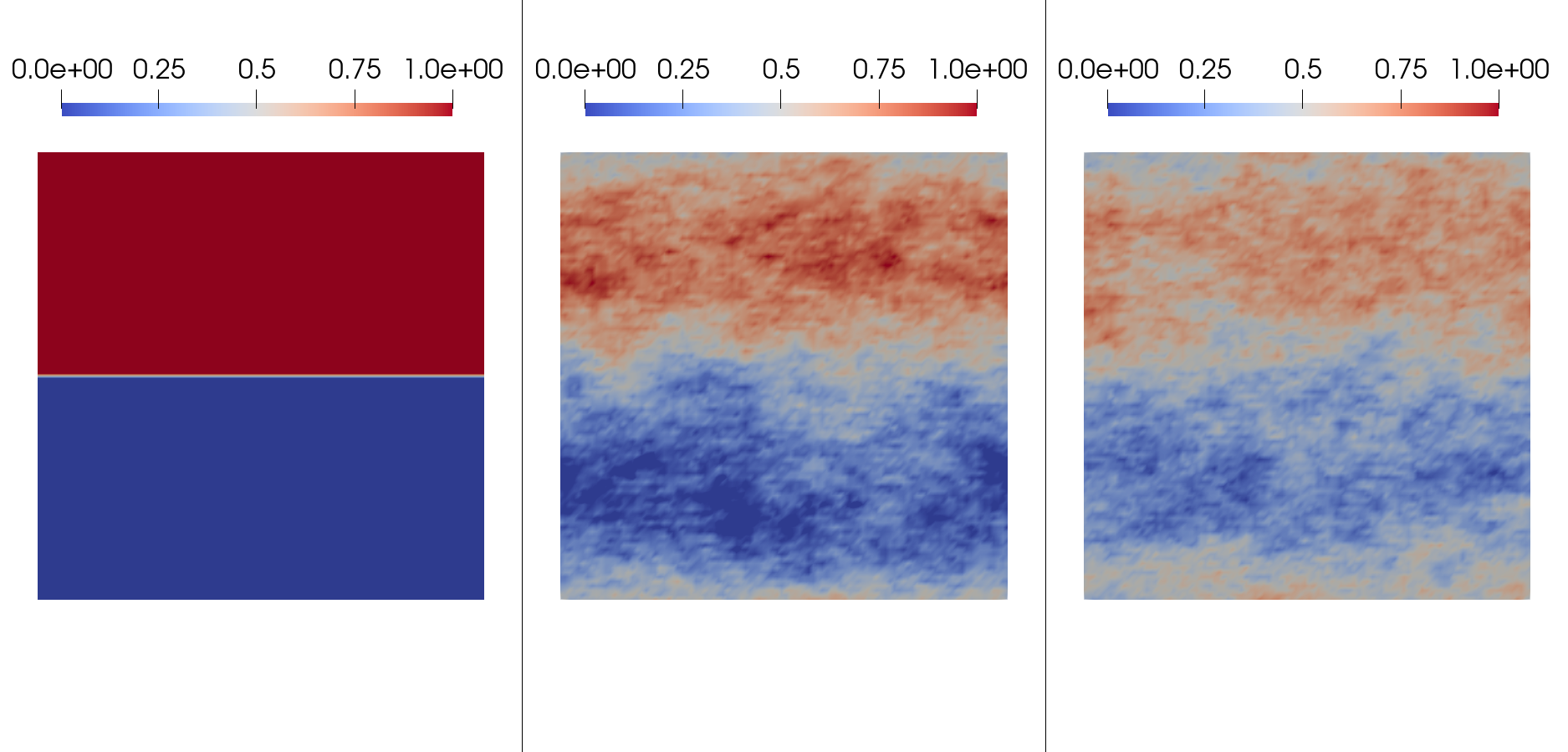}
	\caption{Solution for ${A^I}$ at time $T=0$, $0.1$ and $0.2$}
\end{figure}
We convert the numerical results from $A$ to $u$.
\begin{figure}[H]
	\centering
	\includegraphics[width=12cm]{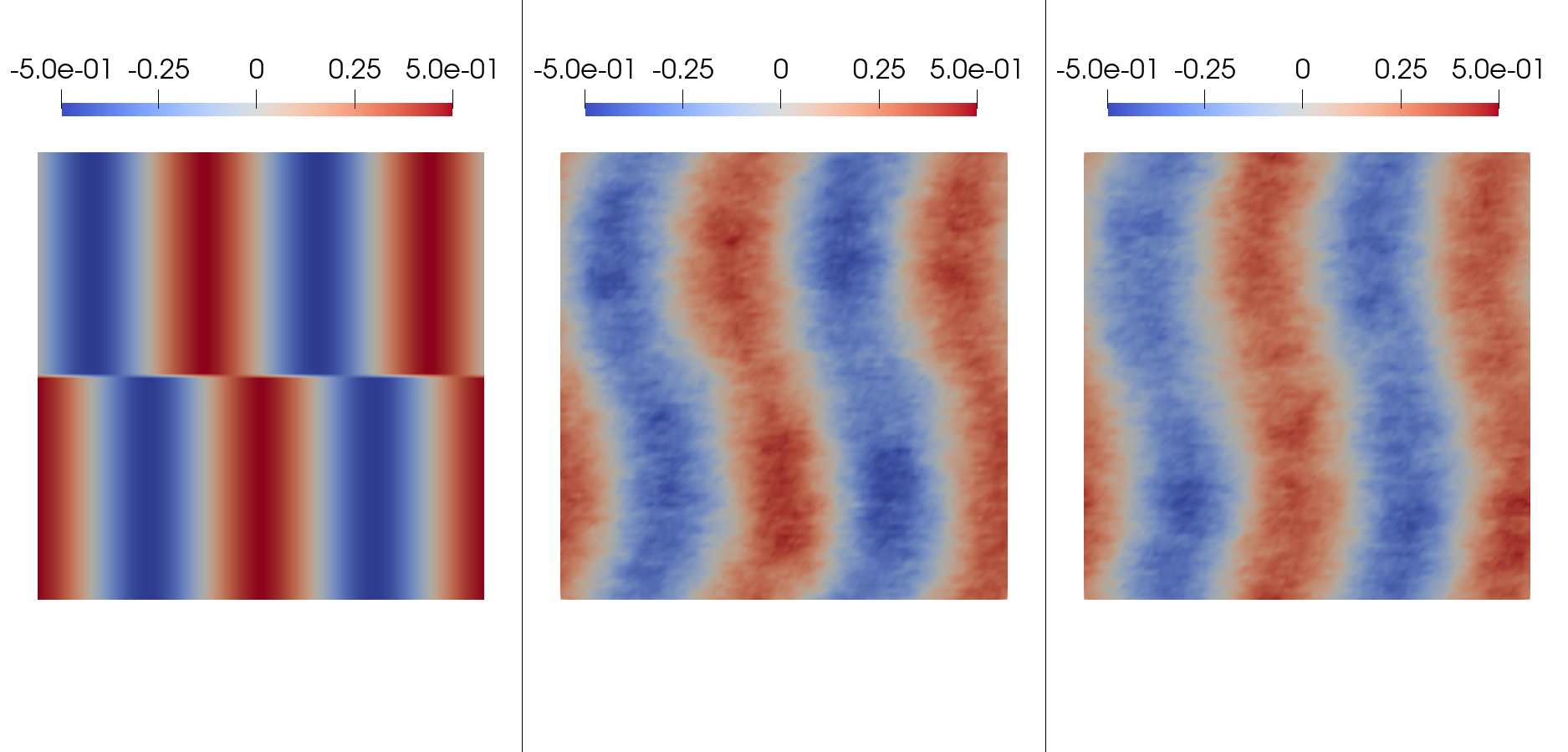}
	\caption{$u$ at time $t=0, 0.1/\varepsilon^2$ and $0.2/\varepsilon^2$ converted by the Ansatz from $A$ of time $T=0$, $0.1$ and $0.2$}
	\label{u_byAnsatz_sto_eps005}
\end{figure}
Then we perform direct simulation for $u$ and compare it to the simulation results of $u$ convert by $A$ which are presented in Figure \ref{u_byAnsatz_sto_eps005}.
\begin{figure}[H]
	\centering
	\includegraphics[width=12cm]{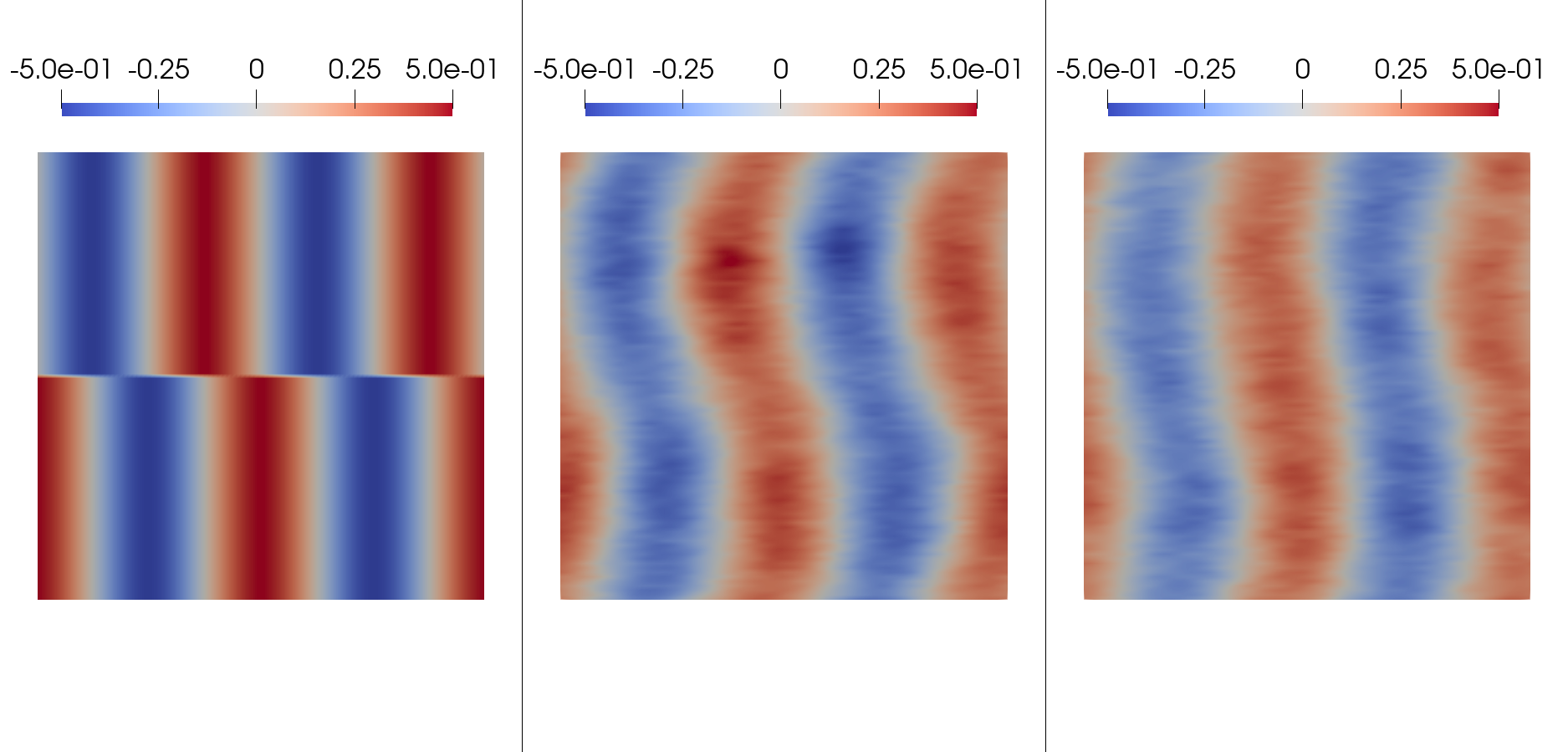}
	\caption{$u$ by direct simulation at time $t=0, 0.1/\varepsilon^2$ and $0.2/\varepsilon^2$}
\end{figure}

\vspace{3mm}

\noindent
{\bf Acknowledgement} This work was supported by JSPS KAKENHI Grant Numbers JP19KK0066, JP20K03669.
The work of Guido Schneider is partially supported by the Deutsche
Forschungsgemeinschaft DFG through the cluster of excellence `SimTech' 
under EXC 2075-390740016.

%
%

\end{document}